\pdfoutput=1
\RequirePackage{ifpdf}
\ifpdf 
\documentclass[pdftex]{sigma}
\else
\documentclass{sigma}
\fi

\usepackage{pictexwd,dcpic}

\DeclareMathOperator{\res}{res} 
\DeclareMathOperator{\TR}{TR} \DeclareMathOperator{\Op}{Op}
\DeclareMathOperator{\Tr}{Tr} \DeclareMathOperator{\ord}{ord}
\DeclareMathOperator{\odd}{odd} \DeclareMathOperator{\supp}{supp}
\DeclareMathOperator{\tr}{tr} \DeclareMathOperator{\Log}{Log}
\DeclareMathOperator{\Lie}{Lie} \DeclareMathOperator{\Det}{Det}
\DeclareMathOperator{\Exp}{Exp} \DeclareMathOperator{\rk}{rk}

\def\dbar{d{\hskip-1pt\bar{}}\hskip1pt}

\newcommand{\R}{{\mathbb R}}
\newcommand{\C}{{\mathbb C}}
\newcommand{\N}{{\mathbb N}}
\newcommand{\Z}{{\mathbb Z}}

\def \Cl {{C\ell}}

\begin{document}

\allowdisplaybreaks

\renewcommand{\PaperNumber}{023}

\FirstPageHeading

\ShortArticleName{Traces and Determinants on Odd-Class Operators}

\ArticleName{Classif\/ication of Traces and Associated Determinants\\ on Odd-Class Operators in Odd Dimensions}

\Author{Carolina NEIRA JIM\'ENEZ~$^\dag$ and Marie Fran{\c{c}}oise OUEDRAOGO~$^\ddag$}

\AuthorNameForHeading{C.~Neira Jim\'enez and M.F.~Ouedraogo}

\Address{$^\dag$~Fakult\"at f\"ur  Mathematik, Universit\"at Regensburg, 93040 Regensburg, Germany}
\EmailD{\href{mailto:Carolina.Neira-Jimenez@mathematik.uni-regensburg.de}{Carolina.Neira-Jimenez@mathematik.uni-regensburg.de}}
\URLaddressD{\url{http://homepages.uni-regensburg.de/~nec07566/}}

\Address{$^\ddag$~D\'epartment de Math\'ematiques, Universit\'e de Ouagadougou, 03 BP 7021, Burkina Faso}
\EmailD{\href{mailto:marie.oued@univ-ouaga.bf}{marie.oued@univ-ouaga.bf}}

\ArticleDates{Received November 30, 2011, in f\/inal form April 11, 2012; Published online April 21, 2012}

\Abstract{To supplement the already known classif\/ication of traces on classical pseudodif\/ferential operators, we present a classif\/ication of traces on the algebras of odd-class pseudodif\/ferential operators of non-positive order acting on smooth functions on a closed odd-dimensional manifold. By means of the one to one correspondence between conti\-nuous traces on Lie algebras and determinants on the associated regular Lie groups, we give a~classif\/ication of determinants on the group associated to the algebra of odd-class pseudo\-dif\/ferential operators with f\/ixed non-positive order. At the end we discuss two possible ways to extend the def\/inition of a determinant outside a neighborhood of the identity on the Lie group associated to the algebra of odd-class pseudodif\/ferential operators of order zero.}


\Classification{58J40; 47C05}

\section{Introduction}

  From the connection between the trace of a matrix with scalar coef\/f\/icients and its eigenvalues, one can derive a relation between the trace and the determinant of a matrix, namely
\begin{gather}
\det(A)=\exp(\tr(\log A)). \label{detexptr}
\end{gather}
At the level of Lie groups, a trace on a Lie algebra is the derivative of the determinant at the identity on the associated Lie group. Using the exponential mapping between a Lie algebra and its Lie group, one recovers in this setting the relation~\eqref{detexptr}. This exponential mapping always exists in the case of f\/inite-dimensional Lie groups, and in the inf\/inite-dimensional case its existence is ensured by requiring regularity of the Lie group.

On trace-class operators over a separable Hilbert space one can promote the trace on matrices to an operator trace. Further generalizing to classical pseudodif\/ferential operators one can consider traces on such operators. In the case of a closed manifold of dimension greater than one, M.~Wodzicki proved that there is a unique trace (up to a constant factor) on the whole algebra of classical pseudodif\/ferential operators acting on smooth functions on the manifold, namely the noncommutative residue \cite{wodzicki}. As S.~Paycha and S.~Rosenberg pointed out~\cite{rosenberg}, this fact does not rule out the existence of other traces when restricting to subalgebras of such operators. In fact, other traces such as the leading symbol trace, the operator trace and the canonical trace appear naturally on appropriate subalgebras. The classif\/ication of the traces on algebras of classical pseudodif\/ferential operators of non-positive order has been carried out in~\cite{neirathesis} (see also~\cite{leschneira}).

After the construction of the canonical trace on non-integer order pseudodif\/ferential operators, M.~Kontsevich and S.~Vishik \cite{kontsevichdeterminant} introduced the set of odd-class operators, which is an algebra that contains the dif\/ferential operators. They also def\/ined a trace on this algebra when the dimension of the manifold is odd, and it has been proven that this is the unique trace in this context \cite{maniccia, sylvie,ponge}. In~\cite{ouedraogothese}, M.F.~Ouedraogo gave another proof of this fact based on the expression of a symbol of such an operator as a sum of derivatives of symbols corresponding to appropriate operators on the same algebra.

Odd-class operators are one of the rare types of operators in odd dimensions on which the canonical trace is def\/ined, this fact serving as a motivation to investigate the classif\/ication of traces on the algebra of odd-class operators of order zero. The noncommutative residue is not of interest here since it vanishes on odd-class operators (Lemma~\ref{eq: nullityRES}). This contrasts with the algebra of ordinary zero-order operators where the only traces are linear combinations of the leading symbol trace and the noncommutative residue (see~\cite{lescure}). We supplement that classif\/ication of traces, showing that when restricting to odd-class zero-order operators, the only traces turn out to be linear combinations of the leading symbol trace and the canonical trace (this is the particular case $a=0$ in Theorem~\ref{classifictracesClodda}).

In this article we present the classif\/ication of traces on algebras of odd-class pseudodif\/ferential operators acting on smooth functions on a closed odd-dimensional manifold. The methods we implement combine various approaches used in the literature on the classif\/ication of traces. However, a detailed analysis is required here because of the specif\/icity of odd-class operators (see Proposition \ref{prop:PDOsumComutators}). We recall the one to one correspondence between continuous traces and determinants of class $C^1$ on regular Lie groups\footnote{Special instances of this one to one correspondence were discussed by P.~de la Harpe and G.~Skandalis in~\cite{HarpeSkandalis}.}, and as in \cite{lescure} we use this correspondence to give the classif\/ication of determinants on the Lie group associated to the algebra of odd-class operators of a f\/ixed non-positive order. At the end we discuss two ways to extend the def\/inition of determinants outside a neighborhood of the identity.

In Section~\ref{section2} we recall some of the basic notions of classical pseudodif\/ferential operators, including that of symbols on an open subset of the Euclidean space and odd-class operators on a closed manifold. Inspired by \cite{maniccia}, we use the representation of an odd-class symbol as a sum of derivatives up to a smoothing symbol (Proposition \ref{prop:SymbolSumDerivatives}), to express an odd-class operator in terms of commutators of odd-class operators (Proposition~\ref{prop:PDOsumComutators}), a fact that helps considerably in the classif\/ication of traces.

In Section~\ref{sectiontraces} we give a classif\/ication of traces on odd-class operators of non-positive order. For that we recall the known traces on classical pseudodif\/ferential operators. The noncommutative residue vanishes on the algebra of odd-class operators in odd dimensions, whereas the canonical trace is the unique linear form on this set which vanishes on commutators of elements in the algebra (see~\cite{maniccia}). Then using the fact that any odd-class operator can be expressed in terms of commutators of odd-class operators (Proposition~\ref{prop:PDOsumComutators}), we prove that any trace on an algebra of odd-class operators of f\/ixed non-positive order can be expressed as a linear combination of a~generalized leading symbol trace and the canonical trace (Theorem~\ref{classifictracesClodda}).

In Section~\ref{section4} we classify determinants on the Lie groups associated to the algebras of odd-class operators of non-positive order. We follow some of the work done in~\cite{lescure}, concerning the one to one correspondence between continuous traces on Lie algebras and $C^1$-determinants on the associated regular Lie groups (this is also discussed in~\cite{HarpeSkandalis} in specif\/ic situations). Then, we combine this correspondence with the classif\/ication of traces given in Section~\ref{section:classificationZero}, to provide the classif\/ication of determinants on Lie groups associated to algebras of odd-class operators of f\/ixed non-positive order (Theorem~\ref{prop:ClassificationDetCl_{odd}^{0}(M)}). This classif\/ication is carried out for operators in a small neighborhood of the identity, where the exponential mapping is a dif\/feomorphism.

 At the end of this section, we give two possible ways to extend the def\/inition of a determinant outside a neighborhood of the identity on the Lie group associated to the algebra of odd-class pseudodif\/ferential operators of order zero; the f\/irst one (see~\eqref{Detlambda}) using a spectral cut to def\/ine the logarithm of an admissible operator; in this case, for some traces this def\/inition of determinant depends on the spectral cut; the other one (see~\eqref{Detlambdagamma}) via the def\/inition of the determinant of an element on the pathwise connected component of the identity, using a path that connects the element with the identity. Both the f\/irst and the second extension of the def\/inition of a determinant, provide maps which do not depend on the spectral cut and which satisfy the multiplicativity property, under the condition that the image of the fundamental group of invertible odd-class pseudodif\/ferential operators of order zero is trivial.

\section{Preliminaries on pseudodif\/ferential operators}\label{section2}

Here we recall the basic notions of classical pseudodif\/ferential operators following \cite{Shubin}.

\subsection{Symbols}

Let $U$ be an open subset of $\mathbb{R}^n$. Given $a\in\mathbb{C}$, a {\it symbol} of order $a$ on $U$ is a complex valued function $\sigma(x,\xi)$ in $C^{\infty}(U\times\mathbb{R}^n)$ such that for any compact subset $K$ of $U$ and any two multiindices $\alpha=(\alpha_1,\dots,\alpha_n), \beta=(\beta_1,\dots,\beta_n) \in \mathbb{N}^n$ there exists a constant $C_{K, \alpha, \beta}$ satisfying for all $(x,\xi)\in K\times\mathbb{R}^n$,
\[
\big| \partial_x^{\alpha}\partial_{\xi}^{\beta}\sigma(x,\xi) \big| \leq C_{K, \alpha, \beta}(1+ \vert \xi \vert )^{\Re(a)- \vert \beta\vert},
\]
where $\partial_x^{\alpha}=\partial_{x_1}^{\alpha_1}\cdots\partial_{x_n}^{\alpha_n}$, $\vert \beta\vert = \beta_1+\cdots+\beta_n$, and $\Re(a)$ stands for the real part of $a$. Let $S^{a}(U)$ denote the set of such symbols.

Notice that if $\Re(a_1)<\Re(a_2)$, then $S^{a_1}(U)\subset S^{a_2}(U)$. We denote by $S^{-\infty}(U):= \bigcap\limits_{a \in \mathbb{C}}S^a(U)$ the space of smoothing symbols on $U$. Given $\sigma\in S^{m_0}(U),\ \sigma_j\in S^{m_j}(U)$, where $m_j\to-\infty$ as $j\to\infty$, we write
\[
\sigma\sim\sum_{j=0}^\infty\sigma_j,
\]
if for every $N\in\mathbb{N}$
\[
\sigma-\sum_{j=0}^{N-1}\sigma_j\in S^{m_N}(U).
\]

The product $\star$ on symbols is def\/ined as follows: if $\sigma_1 \in S^{a_1}(U)$ and $\sigma_2\in S^{a_2}(U)$,
\begin{gather}
\sigma_1\star \sigma_2(x, \xi) \sim \sum_{|\alpha|\geq0}\frac{(-i)^{\vert\alpha\vert}}{\alpha!}\partial_{\xi}^{\alpha} \sigma_1(x,\xi)\partial_x^{\alpha}\sigma_2(x,\xi). \label{eq:starproduct}
\end{gather}
In particular, $\sigma_1\star \sigma_2 \in S^{a_1 + a_2}(U)$.

\subsubsection{Classical symbols}\label{sectionsymbols}

A symbol $ \sigma\in S^a(U)$ is called  {\it classical}, and we write $\sigma\in CS^a(U)$, if there is an asymptotic expansion
\begin{gather}
\sigma(x,\xi) \sim \sum_{j=0}^{\infty}\psi(\xi)\, \sigma_{a-j}(x,\xi). \label{eq:classicallocallogsymb}
\end{gather}
Here $\sigma_{a-j}(x,\xi)$ is a positively homogeneous function in $\xi$ of degree $a-j$:
\[
\sigma_{a-j}(x,t\xi)=t^{a-j}\sigma_{a-j}(x,\xi) \qquad \text{for all}  \ \ t\in\mathbb{R}^+, \ \ \vert \xi\vert \not=0,
\]
and $\psi\in C^{\infty}(\mathbb{R}^n)$ is any cut-of\/f function which vanishes for $\vert \xi \vert \leq \frac{1}{2}$ and such that $\psi(\xi)=1$ for $\vert \xi\vert \geq 1$.

We denote by
\[
CS(U)=\left\langle \bigcup_{a\in \mathbb{C}} CS^a(U)\right\rangle
\]
the algebra generated by all classical symbols on $U$ for the product $\star$.

\subsubsection{Odd-class symbols}

The homogeneous components in the asymptotic expansion of a classical symbol may satisfy some other symmetry relations additional to the positive homogeneity on the second variable. Now we recall the def\/inition of odd-class symbols introduced f\/irst in \cite{kontsevichdeterminant} (see also \cite{grubb}).

\begin{definition}[see \cite{kontsevichdeterminant}]
A classical symbol $\sigma\in {CS}^{a}(U)$ with integer order $a$ is \emph{odd-class} if for each $j\geq 0$, the term $\sigma_{a-j}$ in the asymptotic expansion (\ref{eq:classicallocallogsymb}) satisf\/ies
\begin{gather}
\sigma_{a-j}(x,-\xi)=(-1)^{a-j}\sigma_{a-j}(x,\xi) \qquad {\rm for}\quad \vert \xi\vert \geq 1. \label{eq:OddlogPolyhom}
\end{gather}
In other words, odd-class symbols have the parity one would expect from dif\/ferential operators.
\end{definition}
Let us denote by $CS_{\odd}^a(U)$ the set of odd-class symbols of order $a\in \mathbb{Z}$ on $U$. We set
\[
CS_{\odd}(U)=\bigcup_{a\in \mathbb{Z}} CS_{\odd}^a(U).
\]

\begin{lemma}[see \cite{ducourtiouxthese, kontsevichdeterminant}]\label{lemmaoddclasssymbols}
The odd-class symbols satisfy the following:
\begin{enumerate}\itemsep=0pt
\item The product $\star$ of two odd-class symbols is an odd-class symbol, therefore $CS_{\odd}(U)$ is an algebra.
\item If an odd-class symbol is invertible with respect to the product $\star$, then its inverse is an odd-class symbol.
\end{enumerate}
\end{lemma}

\subsubsection{The noncommutative residue on symbols}

As before, we consider $U$ an open subset of $\R^n$.

\begin{definition}[see \cite{guillemin,wodzicki}]
The \emph{noncommutative residue} of a classical symbol $\sigma\sim\sum\limits_{j=0}^\infty\sigma_{a-j}\in CS^a(U)$ is def\/ined by
\[
\res(\sigma):=\int_U \int_{S_x^{*}U} \sigma_{-n}(x,\xi) \mu(\xi)\, dx=:\int_U \res_x(\sigma)dx,
\]
where $\mu$ is the surface measure on the unit sphere $S_x^{*}U$ over $x$ in the cotangent bundle $T^*U$.
\end{definition}
The noncommutative residue clearly vanishes on symbols of order strictly less than $-n$ and also on symbols of non-integer order.

\begin{lemma}\label{eq: nullityRES}
In odd dimensions, the noncommutative residue of any odd-class symbol vanishes.
\end{lemma}

\begin{proof}
Let $\sigma\in CS_{\odd}^a(U)$ be with asymptotic expansion $\sigma\sim\sum\limits_{j=0}^\infty\psi \sigma_{a-j}$ as in \eqref{eq:classicallocallogsymb}. Since $n$ is odd, we have $\sigma_{-n}(x,-\xi)=(-1)^{n}\sigma_{-n}(x,\xi)=-\sigma_{-n}(x,\xi)$. Then we obtain for any $x\in U$
\begin{gather*}
\res_x(\sigma)=\int_{S_x^{*}U}\!\!\sigma_{-n}(x,\xi)\mu(\xi)=-\int_{S_x^{*}U}\!\!\sigma_{-n}(x,-\xi)\mu(\xi)
=-\int_{S_x^{*}U}\!\!\sigma_{-n}(x,\xi)\mu(\xi)=-\res_x(\sigma).
\end{gather*}
Therefore $\res(\sigma)=0$.
\end{proof}

\subsubsection{An odd-class symbol as a sum of derivatives}

\begin{proposition}[see Lemma~1.3 in \cite{fedosov}, and~\cite{maniccia}]\label{prop:SymbolSumDerivatives}
Let $n\in\Z$ be odd. For any $\sigma\in CS_{\odd}^a(U)$, there exist $\tau_i$ in $CS_{\odd}^{a+1}(U)$ such that
\[
\sigma\sim\sum\limits_{i=1}^n\partial_{\xi_i}\tau_i.
\]
\end{proposition}

\begin{proof}
For a cut-of\/f function $\psi$ as in Section \ref{sectionsymbols} consider
\[
\sigma\sim\sum_{j=0}^\infty\psi \sigma_{a-j},
\]
with $\sigma_{a-j}$ a positively homogeneous function of degree $a-j$ in $\xi$ which satisf\/ies \eqref{eq:OddlogPolyhom}.

 If $a-j\not=-n$, consider the homogeneous function $\tau_{i,a-j+1}:=\frac{\xi_i\sigma_{a-j}(x,\xi)}{a-j+n}$. By Euler's identity we have
\[
\sum_{i=1}^n\partial_{\xi_i}(\tau_{i,a-j+1})(x,\xi)=\sigma_{a-j}(x,\xi).
\]
The homogeneous functions $\tau_{i,a-j+1}$ clearly satisfy \eqref{eq:OddlogPolyhom} for $|\xi|\geq1$:
\begin{gather*}
\tau_{i,a-j+1}(x,t\xi) =t^{a-j+1}\tau_{i,a-j+1}(x,\xi),\qquad \forall\, t>0,\\
\tau_{i,a-j+1}(x,-\xi) =(-1)^{a-j+1}\tau_{i,a-j+1}(x,\xi).
\end{gather*}

 Let $a-j=-n$. In polar coordinates $(r,\omega)\in\mathbb{R}^+\times S^{n-1}$, the Laplacian in $\xi$ reads
\[
\Delta=-\sum_{i=1}^n\partial_{\xi_i}^2=-r^{1-n}\partial_r\big(r^{n-1}\partial_r\big)-r^{-2}\Delta_{S^{n-1}}.
\]
Therefore, for any function $f\in C^\infty(S^{n-1})$,
\[
\Delta\big(f(\omega)r^{2-n}\big)=r^{-n}\Delta_{S^{n-1}}f(\omega).
\]
Since $n$ is odd and $\sigma\in CS_{\odd}^a(U)$, by Lemma \ref{eq: nullityRES} we have $\res(\sigma)=0$. Therefore $\sigma_{-n}(x,\cdot)\upharpoonright_{S^{n-1}}$ is orthogonal to the constants which form the kernel $\ker(\Delta_{S^{n-1}})$. Hence there exists a unique function $h(x,\cdot)\in C^\infty(S^{n-1})$, orthogonal to the constants, such that $\Delta_{S^{n-1}}(h(x,\cdot))=\sigma_{-n}(x,\cdot)\upharpoonright_{S^{n-1}}$. The function $h(x,-\xi)+h(x,\xi)$ is constant and orthogonal to the constants, therefore, $h(x,\cdot)$ is an odd function on $S^{n-1}$.

 We choose a smooth function $\chi$ on $\mathbb{R}$ which vanishes for small $r$ and is equal to 1 for $r\geq1/2$. For $r=|\xi|$, we set
\[
b_{-n}(x,\xi):=\chi(|\xi|)|\xi|^{2-n}h\left(x,\frac{\xi}{|\xi|}\right).
\]
The function $b_{-n}$ is smooth on $U\times\mathbb{R}^n$ and is homogeneous of degree $-n+2$ in $\xi$ for $|\xi|\geq1$. As $\sigma_{-n}(x,\xi)$ vanishes for $x$ outside a compact set, so does $b_{-n}(x,\xi)$. In particular, $b_{-n}$ is a symbol of order $-n+2$ on $U$. Let us def\/ine $\tau_{i,-n+1}:=-\partial_{\xi_i}b_{-n}$. Since $h$ is odd so is $b_{-n}$ and therefore,
\[
\tau_{i,-n+1}(x,-\xi)=-(\partial_{\xi_i}b_{-n})(x,-\xi)
=-\partial_{\xi_i}b_{-n}(x,\xi)
=(-1)^{-n+1}\tau_{i,-n+1}(x,\xi).
\]
Moreover, we have for $r=|\xi|\geq1$
\[
\Delta b_{-n}(x,\cdot)=\Delta\big(r^{2-n}h(x,\cdot)\big)=r^{-n}\left(\sigma_{-n}(x,\cdot)\upharpoonright_{S^{n-1}}\right)=\sigma_{-n}(x,\cdot).
\]

Let $\tau_i\sim\sum\limits_{j=0}^\infty\psi\,\tau_{i,a-j+1}$, then since $\partial_{\xi_i}\psi$ has compact support, the dif\/ference $\sigma-\sum\limits_{i=1}^n\partial_{\xi_i}\tau_i$ is smoothing and
\begin{gather*}
\sigma\sim\sum_{i=1}^n\sum_{j=0}^\infty\psi \partial_{\xi_i}\tau_{i,a-j+1}\sim\sum_{i=1}^n\partial_{\xi_i}\tau_i. \tag*{\qed}
\end{gather*}
\renewcommand{\qed}{}
\end{proof}

\subsection{Pseudodif\/ferential operators}

Let $U\subset \mathbb{R}^n$ be an open subset, and denote by $C^{\infty}_c (U)$ the space of smooth compactly supported functions on $U$. To the symbol $\sigma\in S(U)$, we associate the linear integral operator $\Op(\sigma):C^{\infty}_c (U)\to C^{\infty}(U)$ def\/ined for $u\in C^{\infty}_c (U)$ by
\begin{gather*}
\Op(\sigma)(u)(x) =\int_{T_x^*U}{e^{ix\cdot\xi}\sigma(x, \xi)\widehat u(\xi)\dbar\xi}
 =\int_{T_x^*U}{\int_{U}{e^{i(x-y)\cdot\xi}\sigma(x,\xi) u(y) dy\dbar\xi}}\\
\hphantom{\Op(\sigma)(u)(x)}{} =\int_{U}k(x,y)u(y)dy,
\end{gather*}
where $\widehat u(\xi)=\int_{U}{e^{-iy\cdot\xi}u(y)\,dy}$ is the Fourier transform of $u$ and $\dbar\xi:=(2\pi)^{-n}d\xi$.
In this expression $k(x,y)=\int{e^{i(x-y)\cdot\xi}\sigma(x, \xi)\dbar\xi}$ is seen as a distribution on $U\times U$ that is smooth outside the diagonal. We say that $\Op(\sigma)$ is a {\it pseudodifferential operator} ($\psi$DO) with Schwartz kernel given by $k(x,y)$. An operator is \emph{smoothing} if its Schwartz kernel is a smooth function on $U\times U$. If $\sigma\sim \sum\limits_{j=0}^{\infty}\psi\sigma_{a-j}$ is a classical symbol of order~$a$, then $A=\Op(\sigma)$ is called a {\it classical pseudodifferential operator} of order~$a$. The homogeneous component $\sigma_a$ of $\sigma$ is called {\it the leading symbol} of~$A$, and will be denoted by $\sigma^L_A$.

A $\psi$DO $A$ on $U$ is called {\it properly supported} if for any compact $K\subset U$, the set $\{(x,y)\in \supp(k_A): x\in K  \ {\rm or}\ y\in K\}$ is compact, where $\supp(k_A)$ denotes the support of the Schwartz kernel of~$A$.  Any $\psi$DO $A$ can be written in the form (see~\cite{Shubin})
\begin{gather}
A=P+R, \label{eq:PDO}
\end{gather}
where $P$ is properly supported and $R$ is a smoothing operator.

A properly supported $\psi$DO maps $C^{\infty}_c(U)$ into itself. The product $\star$ on symbols def\/ined in \eqref{eq:starproduct} induces a composition on properly supported $\psi$DOs on $U$. The composition $AB$ of two properly supported $\psi$DOs $A$ and $B$ is a properly supported $\psi$DO with symbol $\sigma(AB)=\sigma(A) \star \sigma(B)$.

The notion of a $\psi$DO can be extended to operators acting on manifolds (see Section~4.3 in~\cite{Shubin}). Let $M$ be a smooth closed manifold of dimension $n$. A linear operator $A:C^{\infty}(M)\to C^{\infty}(M)$ is a \emph{pseudodifferential operator} of order $a$ on $M$ if in any atlas, $A$ is locally a pseudodif\/ferential operator. This means that given a local coordinate chart $U$ of $M$, with dif\/feo\-mor\-phism $\varphi:U\to V$, from $U$ to an open set $V\subseteq\R^n$, the operator $\varphi^\#A$ def\/ined by the following diagram is a pseudodif\/ferential operator of order $a$ on $V$:
\[
\begindc{\commdiag}[35]
\obj(1,3)[objLie(cal{G})]{$C_c^\infty(V)$}
\obj(4,3)[objMathbb{C}]{$C^\infty(V)$}
\obj(4,1)[objMathbb{C}^*]{$C^\infty(U)$}
\obj(1,1)[objCal{G}]{$C_c^\infty(U)$}
\mor{objLie(cal{G})}{objMathbb{C}}{$\varphi^\#A$}
\mor{objMathbb{C}}{objMathbb{C}^*}{$\varphi^*$}
\mor{objCal{G}}{objMathbb{C}^*}{$r_{U}\circ A\circ i_{U}$}
\mor{objLie(cal{G})}{objCal{G}}{$\varphi^*$} [\atright,\solidarrow]
\enddc
\]
where $i_{U}:C_c^\infty(U)\to C_c^\infty(M)$ is the natural embedding, and $r_{U}:C^\infty(M)\to C^\infty(U)$ is the natural restriction.

Let $\Cl^a(M)$ denote the set of classical $\psi$DOs of order $a$ on $M$, i.e.\ operators whose symbol is classical of order $a$ in any local chart of $M$. If $A_1\in \Cl^{a_1}(M), A_2\in \Cl^{a_2}(M)$, then $A_1 A_2\in \Cl^{a_1 + a_2}(M)$, thus, the space $\Cl^{a}(M)$ is an algebra if and only if $a$ is an integer and $a\leq0$. We denote by
\[
\Cl(M):= \left\langle\bigcup_{a\in \mathbb{C}} \Cl^a(M)\right\rangle
\]
the algebra generated by all classical $\psi$DOs on $M$, and by
\[
\Cl^{-\infty}(M):=\bigcap_{a\in \mathbb{C}} \Cl^a(M)
\]
the ideal of smoothing operators in $\Cl(M)$.

We will also denote by
\[
\Cl^{\notin\Z}(M):=\left\langle\bigcup\limits_{a\notin\Z\cap[-n,+\infty)}\Cl^{a}(M)\right\rangle
\]
the space generated by classical $\psi$DOs on $M$ whose order is non-integer or less than $-n$.

A classical operator $A\in \Cl^{a}(M)$ of integer order $a$ is odd-class if in any local chart its local symbol $\sigma(A)$ is odd-class.
We denote by $\Cl_{\odd}^a(M)$ the set of odd-class operators of order $a\in \mathbb{Z}$ and we def\/ine
\[
\Cl_{\odd}(M)=\bigcup_{a\in \mathbb{Z}} \Cl_{\odd}^a(M).
\]

As in Lemma \ref{lemmaoddclasssymbols}, the following lemma implies that $\Cl_{\odd}(M)$ is an algebra:
\begin{lemma}[see Section~4 in~\cite{kontsevichdeterminant}]\label{lema:oddClassproduct}
Let $A\in \Cl^{a}_{\odd}(M)$ and $B\in \Cl^{b}_{\odd}(M)$, $a,b\in\Z$. Then $AB\in \Cl^{a+b}_{\odd}(M)$. If moreover $B$ is an invertible elliptic operator, then $B^{-1}\in \Cl^{-b}_{\odd}(M)$ and $AB^{-1}\in \Cl^{a-b}_{\odd}(M)$.
\end{lemma}

The algebra $\Cl_{\odd}(M)$ contains the dif\/ferential operators and their parametrices.

\begin{remark}\label{remarknotationClodd}
Even though the def\/inition of odd-class pseudodif\/ferential operators makes sense on any closed manifold, in this paper we restrict ourselves to odd-dimensional closed manifolds. The reason is that the canonical trace (which will be explained below in Section~\ref{subsectioncanonicaltrace}) is well def\/ined only in that case. So, from now on, the notation $\Cl_{\odd}(M)$ will be used only when the dimension $n$ of the manifold $M$ is odd.
\end{remark}

\subsubsection{Fr\'echet topology on pseudodif\/ferential operators}\label{subsectionfrechettop}

For any complex number $a$,  we equip the vector space $\Cl^a(M)$ with a Fr\'echet topology as follows. Let us consider a covering of $M$ by open neighborhoods $\{U_i\}_{i\in I}$, a f\/inite subordinated partition of unity $(\chi_i)_{i\in I}$ and smooth functions $(\widetilde{\chi_i})_{i\in I}$ on $M$ such that $\supp(\widetilde{\chi_i})\subset U_i$ and $\widetilde{\chi_i}=1$ near the support of $\chi_i$. By \eqref{eq:PDO} any $\psi$DO $A$ can be written in the form $A=\sum\limits_{i\in I}(A_i+R_i)$ where the operators $A_i:=\chi_i\cdot\Op(\sigma_i)\cdot\widetilde{\chi_i}\in Cl^{a}(M)$ are properly supported in $U_i$ with symbols $\sigma^{(i)}(A)\sim\sum\limits_{j=0}^\infty\psi \sigma_{a-j}^{(i)}(A)$, and $R_i$ is a smoothing operator with smooth kernel $k_i$ which has compact support in $U_i\times U_i$.

We equip ${\Cl}^{a}(M)$ with the following countable set of semi-norms: for any compact subset $K\subset U_i$ for any $j\geq 0$, $N\geq 1$ and for any multiindices $\alpha$, $\beta$
\begin{gather*}
 \sup\limits_{x\in K}\sup\limits_{
  \xi\in\mathbb{R}^n}(1+\vert\xi\vert)^{\vert\beta\vert-a}
  \big|\partial_x^{\alpha}\partial_{\xi}^{\beta} \sigma^{(i)}(A)(x,\xi)\big|,\\
 \sup\limits_{x\in K}\sup\limits_{
  \xi\in\mathbb{R}^n}(1+\vert\xi\vert)^{\vert\beta\vert-a+N}
  \left|\partial_x^{\alpha}\partial_{\xi}^{\beta} \left(\sigma^{(i)}(A)-\sum_{j=0}^{N-1}\chi(\xi)  \sigma_{a-j}^{(i)}(A)\right)(x,\xi)\right|,\\
 \sup\limits_{x\in K}\sup\limits_{ \vert\xi\vert=1}
  \big|\partial_x^{\alpha}\partial_{\xi}^{\beta}\sigma_{a-j}^{(i)}(A)(x,\xi)\big|,\qquad
 \sup\limits_{x,y\in K}
  \big|\partial_x^{\alpha}\partial_{y}^{\beta} k_i(x,y)\big|.
\end{gather*}

\subsubsection{The logarithm of a classical pseudodif\/ferential operator}

An operator $A \in \Cl(M)$ with positive order has \emph{principal angle} $\theta$ if for every $(x, \xi)\in T^*M\setminus\{0\}$, the leading symbol $\sigma_A^L(x, \xi)$ has no eigenvalues on the ray $L_\theta=\{re^{i\theta}, r\geq 0\}$; in that case $A$ is \emph{elliptic}.
\begin{definition}[see e.g.~\cite{okikiolu}]
An operator $A\in \Cl(M)$ is \emph{admissible} with spectral cut $\theta$ if $A$ has principal angle $\theta$ and the spectrum of $A$ does not meet $L_{\theta}=\{re^{i\theta}, r\geq 0\}$. In particular such an operator is invertible and elliptic. The angle $\theta$ is called an \emph{Agmon angle} of $A$.
\end{definition}
Let $A\in \Cl(M)$ be admissible with spectral cut $\theta$ and positive order $a$. For $\Re(z)<0$, the complex power $A_{\theta}^z$ of $A$ is
def\/ined by the Cauchy integral (see \cite{seeley})
\[
A_\theta^{z}=\frac{i}{2\pi} \int_{\Gamma_{r,\theta}} \lambda_\theta^{z} (A-\lambda)^{-1}\, d\lambda,
\]
where $\lambda_\theta^z= \vert \lambda\vert^z e^{iz {(\rm arg}\lambda)}$ with  $\theta\leq {\rm arg}\,\lambda<\theta+2\pi$. Here $\Gamma_{r,\theta}$ is a contour along the ray $L_\theta$ around the (non-zero) spectrum of $A$, and $r$ is any small positive real number such that $\Gamma_{r,\theta}$ does not meet the spectrum of $A$. The operator $A_{\theta}^z$ is an elliptic classical $\psi$DO of order $az$; in particular, for $z=0$, we have $ A_{\theta}^0 = I$.

The def\/inition of complex powers can be extended to the whole complex plane by setting $A_{\theta}^z:=A^kA_{\theta}^{z-k}$ for $k\in \N$ and $\Re(z)<k$; this def\/inition is independent of the choice of $k$ in $\N$ and preserves the usual properties, i.e.\ $A_{\theta}^{z_1}A_{\theta}^{z_2}=A_{\theta}^{z_1+z_2}$, $A_{\theta}^k=A^k$, for $k\in\Z$.
Complex powers of operators depend on the choice of spectral cut. Indeed, if $L_{\theta}$ and $L_{\phi}$ are two spectral cuts for~$A$ outside an angle which contains the spectrum of $\sigma^L(A)(x,\xi)$ then
\begin{gather}
A_\theta^z-A_\phi^z=\big(1-e^{2i\pi z} \big)  P_{\theta, \phi} (A) A_\theta^z, \label{eq:DifferenceComplexPowers}
\end{gather}
where the operator
\begin{gather}
 P_{\theta, \phi} (A)=\frac{1}{2i\pi} \int_{\Gamma_{\theta,\phi}}\lambda^{-1}A(A-\lambda)^{-1}\,d\lambda \label{projection}
\end{gather}
is a projection (see \cite{PongeAsymmetry, WodzickiThesis}). Here $\Gamma_{\theta,\phi}$ is a contour separating the part of the spectrum of $A$ contained in the open sector $\theta<\arg\lambda<\phi$ from the rest of the spectrum.

The logarithm  of an admissible operator $A$ with spectral cut $\theta$  is def\/ined in terms of the derivative at $z=0$ of its complex power:
\[
\log_\theta(A)= \partial_z {A_\theta^{z}} _{\vert_{z=0}}.
\]
Logarithms of classical $\psi$DOs of positive order are not classical anymore since their symbols involve a logarithmic term $\log\vert \xi\vert$ as the following elementary result shows.

\begin{proposition}[\cite{okikiolu}]\label{prop:leadsymbolsigmathetaA}
Let $A\in \Cl^{a}(M, E)$ be an admissible operator with positive order and spectral cut $\theta$. In a local trivialization, the symbol of $\log_\theta(A)$ reads:
\[
\sigma_{\log_\theta(A)}(x, \xi)=a  \log \vert \xi\vert  I +\sigma_0^A(x,\xi),
\]
where $\sigma_0^A$ is a classical symbol of order zero.
\end{proposition}

\begin{remark}
If $A$ is a classical $\psi$DO of order zero then $A$ is bounded, and if it admits a spectral cut, then complex powers and the logarithm of~$A$ are directly def\/ined using a Cauchy integral formula, and they are classical $\psi$DOs (see~\cite{okikiolu} and Remark~2.1.7 in~\cite{ouedraogothese}).
\end{remark}

Just as complex powers, the logarithm depends on the choice of spectral cut. Indeed, given two spectral cuts $\theta,\, \phi$ of the operator $A$ such that $0\leq \theta <\phi <2\pi$, dif\/ferentiation of (\ref{eq:DifferenceComplexPowers}) with respect to $z$ and evaluation at $z=0$ yields
\begin{gather}
\log_\theta A-\log_\phi A=-2i\pi P_{\theta, \phi}(A). \label{eq:difflogarithms}
\end{gather}

\subsubsection{Pseudodif\/ferential operators in terms of commutators}

In this subsection we use the $\psi$DO analysis techniques similar to the ones implemented in \cite{maniccia}; we assume that $M$ is an $n$-dimensional closed manifold and $n$ is odd. Given a function $f\in C^\infty(M)$, we also denote by $f$ the zero-order classical $\psi$DO given by multiplication by $f$.
\begin{proposition}\label{prop:PDOsumComutators}
If $A\in\Cl^{a}_{\odd}(M)$, then there exist functions $\alpha_k\in C^\infty(M)$, operators $B_k$ in $\Cl^{a+1}_{\odd}(M)$ and a smoothing operator $R_A$ such that
\begin{gather}
A=\sum_{k=1}^N[\alpha_k,B_k]+R_A. \label{eq:PDOsumComutators}
\end{gather}
\end{proposition}

\begin{proof}
Let us consider a covering of $M$ by open neighborhoods $\{U_j\}_{j=1}^N$ and a f\/inite subordinated partition of unity $\{\varphi_j\}_{j=1}^N$, such that for every pair $(j,k)$, both $\varphi_j$ and $\varphi_k$ have support in one coordinate neighborhood. We write (see Subsection \ref{subsectionfrechettop})
\[
A=\sum_{j,k}\varphi_jA\varphi_k+R.
\]
Each operator $\varphi_jA\varphi_k$ may be considered as an odd-class $\psi$DO on $\mathbb{R}^n$ with symbol $\sigma$ in $CS^a_{\odd}(\mathbb{R}^n).\!$ By Proposition \ref{prop:SymbolSumDerivatives}, there exist odd-class symbols $\tau_l$ of order $a+1$ such that
\[
\sigma\sim\sum\limits_{l=1}^n\partial_{\xi_l}\tau_l.
\]
For any symbol $\tau$ we have,
\[
\sigma([x_l,\Op(\tau)])\sim x_l\cdot \tau-\tau\cdot x_l-i^{-1}\partial_{\xi_l}\tau =i\partial_{\xi_l}\tau,
\]
so that
\[
\Op(\partial_{\xi_l}\tau)=-i[x_l,\Op(\tau)] \quad {\rm up\ to\ a\ smoothing\ operator.}
\]
Since $\sigma\sim\sum\limits_{l=1}^n\partial_{\xi_l}\tau_l,$ there exist smoothing operators $R''$, $R^\prime$ such that
\[
\Op(\sigma)=\Op\left(\sum\limits_{l=1}^n\partial_{\xi_l}\tau_l\right)+R''=-i\sum\limits_{l=1}^n [x_l,\Op(\tau_l)]+R^\prime.
\]
For each index $j$ let $\psi_j\in C_c^\infty(U_j)$ be such that $\psi_j\equiv1$ near $\supp(\varphi_j)$. Then for $\chi$ in $C_c^\infty(\mathbb{R}^n)$ such that $\chi\varphi_j=\varphi_j$, $\chi\psi_j=\psi_j$, we have
\begin{gather*}
\varphi_j\Op(\sigma)\psi_j
 =-i\sum\limits_{l=1}^n\varphi_j[x_l,\Op(\tau_l)]\psi_j+\varphi_jR^\prime\psi_j
 =-i\sum\limits_{l=1}^n [\chi x_l,\varphi_j\Op(\tau_l)\psi_j]+\varphi_jR^\prime\psi_j.
\end{gather*}
As in (\ref{eq:PDO}) we write $\varphi_j A \psi_j= \Op(\sigma_{j}) +R_{j}$ for some $\sigma_j\in CS^{a}_{\odd}(\mathbb{R}^n)$ and some smoothing ope\-ra\-tor~$R_{j}$. We have $A= \sum _{j}\Op(\sigma_{j})+\sum _{j}R_{j}+R$, hence $A$ can be written in the form
\[
A=\sum\limits_{k=1}^{N}[\alpha_k,B_k]+R_A,
\]
where $\alpha_k$ is a smooth function on $M$ (and represents the operator in $\Cl^{0}_{\odd}(M)$ of multiplication by $\alpha_k$), $B_k$ lies in $\Cl^{a+1}_{\odd}(M),$ and $R_A$ is a smoothing operator.
\end{proof}

\begin{corollary}\label{corol2}
If $A\in \Cl^{a}_{\odd}(M)$ then there exist $B_i\in \Cl^{a}_{\odd}(M)$, smooth functions $a_i\in C^\infty(M)$, and a smoothing operator $R$ such that
\[
A-\Op\big(\sigma_A^L\big)=\sum_{i=1}^n[a_i,B_i]+R.
\]
\end{corollary}

\begin{proof}
It follows by applying Proposition~\ref{prop:PDOsumComutators} to $A-\Op(\sigma_A^L)\in \Cl^{a-1}_{\odd}(M)$.
\end{proof}

\section{Classif\/ication of traces on odd-class operators}\label{sectiontraces}

The classif\/ication made in this section is essentially based on Proposition \ref{prop:PDOsumComutators}, namely the decomposition of an odd-class operator in terms of commutators of odd-class operators. Let us f\/irst recall the def\/inition of a trace.

Let $M$ be a closed connected manifold of dimension $n>1$, and let $\mathcal{A}\subseteq\Cl(M)$. A~\emph{trace} on~$\mathcal{A}$ is a map
\[
\tau: \ \mathcal{A}\rightarrow\mathbb{C},
\]
linear in the sense that for all $a, b\in\C$, whenever $A$, $B$ and $aA+bB$ belong to $\mathcal{A}$ we have
\[
\tau(aA+bB)=a \tau(A)+b \tau(B),
\]
and such that for any $A,B\in\mathcal{A}$, whenever $AB,BA\in\mathcal{A}$ it satisf\/ies
\[
\tau([A,B])=0,\qquad \text{or equivalently,} \qquad \tau(AB)=\tau(BA).
\]

\subsection{Examples of traces on pseudodif\/ferential operators}

Interestingly by Lemma \ref{eq: nullityRES}, the noncommutative residue, which is the only trace on the whole space $\Cl(M)$ (see \cite{fedosov,lesch,wodzicki}), vanishes on odd-class operators when the dimension of the manifold is odd. In this section we review the traces which are non-trivial on this class of operators.

\subsubsection[The $L^2$-trace]{The $\boldsymbol{L^2}$-trace}\label{sectionL2trace}

A $\psi$DO $A$ whose order has real part less than $-n$ is a trace-class operator. The $L^2$-trace (also called operator trace) is the functional
\begin{gather}
\Tr:\ \left\langle\bigcup\limits_{\Re(a)<-n}\Cl^{a}(M)\right\rangle \to\mathbb{C},\qquad
A \mapsto\Tr(A):=\int_{M}k_A(x,x)\,dx, \label{TrL2kernel}
\end{gather}
where $k_A$ is the Schwartz kernel of the operator $A$. This trace is continuous for the Fr\'echet topology on the space of $\psi$DOs of constant order less than $-n$.

This is the unique trace on the algebra of smoothing operators $\Cl^{-\infty}(M)$, since we have the exact sequence (see~\cite{guilleminfourier})
\[
0\to\big[\Cl^{-\infty}(M),\Cl^{-\infty}(M)\big]\to \Cl^{-\infty}(M)\overset{\Tr}{\longrightarrow}\mathbb{C}\to0.
\]
More precisely we have

\begin{theorem}[Theorem~A.1 in~\cite{guilleminfourier}]\label{smooth}
If $R$ is a smoothing operator then, for any pseudodifferential idempotent $J$, of rank~$1$, there exist smoothing operators $S_1,\ldots,S_{N}$, $T_1,\ldots,T_{N}$, such that
\[
R=\Tr(R)J+\sum_{j=1}^{N}[S_j,T_j].
\]
Therefore, any smoothing operator with vanishing $L^2$-trace is a sum of commutators in the space $[\Cl^{-\infty}(M),\Cl^{-\infty}(M)]$.
\end{theorem}

\begin{proposition}[see e.g.\ Introduction in~\cite{lesch} and Proposition~4.4 in~\cite{leschtraces}]\label{nocont}
The trace~$\Tr$ does not extend to a trace functional neither on the whole algebra $\Cl(M)$, nor does it on the algebra~$\Cl^0(M)$.
\end{proposition}

\subsubsection{The canonical trace}\label{subsectioncanonicaltrace}

We start with the def\/inition of the cut-of\/f integral of a symbol, as in \cite{sylvie}. Let $U$ be an open subset of $\mathbb{R}^n$.

\begin{proposition}[see \cite{sylvie} and Section~1.2 in~\cite{scott}]
Let $\sigma\in CS^a(U)$, such that for any $N\in\mathbb{N}$, $\sigma$ can be written in the form $ \sigma(x,\xi) =\sum\limits_{j=0}^{N-1}\psi(\xi)\sigma_{a-j}(x,\xi)+\sigma_N(x,\xi)$, with $\sigma_{a-j}$, $\psi$ and $\sigma_N\in S^{a-N}(U)$ as in~\eqref{eq:classicallocallogsymb}. The integral of $\sigma$ over $B^{*}_x(0,R)$ $($which stands for the ball of radius $R$ in the cotangent space $T^{*}_xU)$ has the asymptotic expansion
\[
\int_{B^{*}_x(0,R)}\sigma(x,\xi)d\xi \underset{R\to\infty}{\sim}\sum_{\substack{j=0\\ a-j+n\not=0}}^{N-1}\alpha_j(\sigma)(x)R^{a-j+n} +\res_x(\sigma)\log R+\alpha_x(\sigma),
\]
where $\alpha_x(\sigma)$ converges when $R\to\infty$.
\end{proposition}

\begin{proof}
We write
\[
\int_{B_x^*(0,R)}\psi(\xi)\sigma_{a-j}(x,\xi)d\xi=\int_{B_x^*(0,1)}\psi(\xi)\sigma_{a-j}(x,\xi)d\xi +\int_{B_x^*(0,R)\setminus B_x^*(0,1)}\psi(\xi)\sigma_{a-j}(x,\xi)d\xi.
\]
Using the fact that $\sigma_{a-j}$ is positively homogeneous of degree $a-j$ we have
\begin{gather*}
 \int_{B_x^*(0,R)\setminus B_x^*(0,1)}\psi(\xi)\sigma_{a-j}(x,\xi)d\xi=\int_1^R\int_{\vert\omega\vert=1}r^{a-j+n-1}\sigma_{a-j}(x,\omega)\,d\omega dr \\
\qquad{} =\dfrac{1}{a-j+n}R^{a-j+n} \int_{\vert\omega\vert=1}\sigma_{a-j}(x,\omega)d\omega-\dfrac{1}{a-j+n}\int_{\vert\omega\vert=1}\sigma_{a-j}(x,\omega)\,d\omega.
\end{gather*}
It follows that the integral $\int_{B_x^*(0,R)}\sigma(x,\xi)\,d\xi$ admits the asymptotic expansion
\begin{gather*}
\int_{B_x^*(0,R)}\sigma(x,\xi)\,d\xi\underset{R\to\infty}{\sim}\sum_{\substack{j=0\\a-j+n\not=0}}^{N-1} \frac{1}{a-j+n}R^{a-j+n}\int_{\vert\omega\vert=1}\sigma_{a-j}(x,\omega)\,d\omega \\
\hphantom{\int_{B_x^*(0,R)}\sigma(x,\xi)\,d\xi\underset{R\to\infty}{\sim}}{}
+{\res}_x(\sigma)\log R + \alpha_x(\sigma),
\end{gather*}
where $\alpha_x(\sigma)$ is given by
\begin{gather}
\alpha_x(\sigma) := \int_{B_x^*(0,R)}\sigma_N(x, \xi)\,d\xi+\sum_{j=0}^{N-1} \int_{B_x^*(0,1)} \psi(\xi)\sigma_{a-j}(x,\xi)\,d\xi\nonumber\\
\hphantom{\alpha_x(\sigma) :=}{}
-\sum_{\substack{j=0\\ a-j+n\not=0}}^{N-1} \dfrac{1}{a-j+n}\int_{\vert\omega\vert=1}\sigma_{a-j}(x,\omega)\,d\omega,\label{eq:CutoffIntegral}
\end{gather}
which may depend on the variable $x$.
\end{proof}

 For $N$ suf\/f\/iciently large, the integral $\int_{B_x^*(0,R)}\sigma_N(x,\xi)d\xi$ is well def\/ined, so the term $\alpha_x(\sigma)$ given by~\eqref{eq:CutoffIntegral} converges when $R\to\infty$, and taking this limit we def\/ine the \emph{cut-off integral} of $\sigma$ by
\begin{gather*}
{-\hskip -11,5pt\int}_{T^{*}_xU}\sigma(x,\xi)\,d\xi :=\int_{T^{*}_xU}\sigma_N(x,\xi)\,d\xi+\sum_{j=0}^{N-1} \int_{B^{*}_x(0,1)}\psi(\xi)\sigma_{a-j}(x,\xi)\,d\xi\\
\hphantom{{-\hskip -11,5pt\int}_{T^{*}_xU}\sigma(x,\xi)\,d\xi :=}{}  -\sum_{\substack{j=0\\ a-j+n\not=0}}^{N-1}\dfrac{1}{a-j+n}\int_{S^{*}_xU}\sigma_{a-j}(x,\omega)\,d\omega,
\end{gather*}
where $S^{*}_xU$ stands for the unit sphere in the cotangent space $T^{*}_xU$. This def\/inition is independent of $N>a+n-1$. Moreover, if $a<-n$, then ${-\hskip -10pt\int}_{\mathbb{R}^n}\sigma(x,\xi)d\xi=\int_{\mathbb{R}^n}\sigma(x,\xi)d\xi$.

According to Proposition \ref{nocont}, there is no non-trivial trace on $\Cl(M)$ which extends the $L^2$-trace. However, the $L^2$-trace does extend to non-integer order operators and to odd-class operators. Indeed, M.~Kontsevich and S.~Vishik~\cite{kontsevichdeterminant} constructed such an extension, the \emph{canonical trace}
\begin{gather*}
\TR: \ \Cl^{\notin\Z}(M)\bigcup\Cl_{\odd}(M) \to\mathbb{C},\qquad
 A  \mapsto\TR(A):=\dfrac{1}{(2\pi)^n}\int_Mdx\ {-\hskip -11,5pt\int}_{T_x^{*}M}\sigma(A)(x,\xi)d\xi,
\end{gather*}
where the right hand side is def\/ined using a f\/inite covering of $M$, a partition of unity subordinated to it and the local representation of the symbol, but this def\/inition is independent of such choices. As we already stated in Remark~\ref{remarknotationClodd}, the canonical trace is well def\/ined on $\Cl_{\odd}(M)$ only when the dimension $n$ of the manifold is odd (see~\cite{grubb,lesch}), which is always our case.

 If $A\in \Cl^{a}(M)$, $B\in \Cl^b(M)$ and if $a,b\notin\Z$, then $\ord(AB)=a+b$ may be an integer, so the linear space $\Cl^{\notin\Z}(M)$ is not an algebra; in spite of this, the canonical trace has the following properties (see~\cite{kontsevichdeterminant}, Section~5 in~\cite{lesch}, and \cite{ouedraogo,sylvie,scott}):
 \begin{enumerate}\itemsep=0pt
\item For any $c\in\mathbb{C}$, if $A,B\in\Cl^{\notin\Z}(M)$ are such that $\ord(c A+B)\notin\Z\cap[-n,+\infty)$, or if $A,B\in\Cl_{\odd}(M)$, then $\TR(c A+B)=c\TR(A)+\TR(B)$.
\item For any $A\in\Cl(M)$ such that $\ord(A)<-n$, $\TR(A)=\Tr(A)$, i.e.\ the canonical trace extends the $L^2$-trace def\/ined in \eqref{TrL2kernel}.
\item If $A,B\in\Cl^{\notin\Z}(M)$ are such that $AB\in\Cl^{\notin\Z}(M)$, or if $A,B\in\Cl_{\odd}(M)$, then $\TR(AB)=\TR(BA)$.
\end{enumerate}

The canonical trace is continuous for the Fr\'echet topology on the sets of $\psi$DOs of constant order in $\Cl^{\notin\Z}(M)\bigcup\Cl_{\odd}(M)$.

\subsubsection{Generalized leading symbol traces}\label{subsectiongeneralizedleadingsymboltrace}

In \cite{rosenberg}, S.~Paycha and S.~Rosenberg introduced the leading symbol traces def\/ined on an algebra of operators $\Cl^{a}(M)$ for $a\leq0$; in this section we follow \cite{neirathesis} and consider a trace which actually coincides with a leading symbol trace when $a=0$. Let $a$ be a non-positive integer and consider the projection map $\pi_a$ from $\Cl^{a}(M)$ to the quotient space $\Cl^{a}(M)/\Cl^{2a-1}(M)$:
\begin{gather}
0\to \Cl^{2a-1}(M)\to \Cl^{a}(M)\overset{\pi_a}{\rightarrow}\Cl^{a}(M)/\Cl^{2a-1}(M)\to0. \label{ses1}
\end{gather}
One can see in Remark~4.3.2 of~\cite{neirathesis} that it is possible to construct a splitting
\[
\theta_a: \ \Cl^{a}(M)/\Cl^{2a-1}(M)\to \Cl^{a}(M).
\]

\begin{lemma}
Any continuous linear map $\rho$ on $\Cl^{a}(M)/\Cl^{2a-1}(M)$ defines a trace on $\Cl^{a}(M)$ by $\rho\circ\pi_a$ called \emph{generalized leading symbol trace}.
\end{lemma}

\begin{proof}
If $A, B\in \Cl^{a}(M)$, their commutator $[A,B]$ belongs to $\Cl^{2a-1}(M)$, and since $\rho\circ\pi_a$ vanishes on $\Cl^{2a-1}(M)$, it def\/ines a trace on $\Cl^{a}(M)$.\end{proof}

 For $A\in \Cl^{a}(M)$, $\rho(\pi_a(A))$ depends on $\sigma_a(A),\ldots,\sigma_{2a}(A)$, where $\sigma_{a-i}(A)$ represents the homogeneous component of degree $a-i$ in the asymptotic expansion of the symbol of $A$. Since $\rho\circ\pi_a$ is linear in $A$, it is a linear combination of linear maps $\rho_{a-i}$ on $C^\infty(S^*M)$, in the terms $\sigma_{a-i}(A)$, hence it reads,
\[
\rho(\pi_a(A))=\sum_{i=0}^{|a|}\rho_{a-i}(\sigma_{a-i}(A)).
\]

\begin{remark}
For $a<0$, a leading symbol trace is the particular case of a generalized leading symbol trace when $\rho_{a-i}\equiv0$ for all $i=1,\ldots,|a|$.
\end{remark}

Generalized leading symbol traces are continuous for the Fr\'echet topology on the space of constant order $\psi$DOs, since they are def\/ined in terms of a f\/inite number of homogeneous components of the symbols of the operators.

\subsection[Trace on $\Cl_{\odd}(M)$]{Trace on $\boldsymbol{\Cl_{\odd}(M)}$}

The canonical trace is the unique trace (up to a constant) on its domain. This result was proved in~\cite{maniccia} (which goes back to 2007 but it was published in 2008) using the fact that the operator corresponding to the derivative of a symbol is, up to a smoothing operator, proportional to the commutator of appropriate operators. The latter idea was considered in~\cite{sylvie} to show the equivalence between Stokes' property for linear forms on symbols and the vanishing of linear forms on commutators of operators. In~\cite{ponge} the author uses the Schwartz kernel representation of an operator to express any non-integer order operator and any operator of regular parity class as a sum of commutators up to a smoothing operator, and then to give another proof of the uniqueness of the canonical trace. The following proof, that we f\/ind in Proposition~3.2.4 of~\cite{ouedraogothese}, is done in the spirit of~\cite{maniccia}.

\begin{theorem}\label{prop:UniquenessTrace}
Any trace on $\Cl_{\odd}(M)$ is proportional to the canonical trace.
\end{theorem}

\begin{proof}
By \eqref{eq:PDOsumComutators} any operator $A$ in $\Cl_{\odd}(M)$ can be written in the form
\[
A=\sum\limits_{k=1}^{N}[\alpha_k,B_k]+R_A,
\]
where $\alpha_k$ are smooth functions on $M$ that can be seen as elements of $\Cl^0_{\odd}(M)$, $B_k$ belong to $\Cl^{a+1}_{\odd}(M),$ and $R_A$ is a smoothing operator.
By Theorem \ref{smooth}, we can express $R_A$ as
\[
R_A=\Tr(R_A)J+\sum_{j=1}^{N'}[S_j,T_j],
\]
where $J$ is a pseudodif\/ferential projection of rank 1 and $S_j$, $T_j$ are smoothing operators. Summing up, the expression for $A$ becomes
\begin{gather}
A=\sum\limits_{k=1}^{N}[\alpha_k,B_k]+\Tr(R_A)J+\sum_{j=1}^{N'}[S_j,T_j]. \label{commutators}
\end{gather}
Applying the canonical trace $\TR$ to both sides of this expression, since $\TR$ vanishes on commutators of operators in $\Cl_{\odd}(M)$, we infer that
\[
\TR(A)=\Tr(R_A).
\]
Thus \eqref{commutators} reads
\begin{gather}
A=\sum\limits_{k=1}^{N}[\alpha_k,B_k]+\TR(A)J+\sum_{j=1}^{N'}[S_j,T_j]. \label{commutators2}
\end{gather}
If $\tau$ is a trace on $\Cl_{\odd}(M)$, applying $\tau$ to both sides of~(\ref{commutators2}) we reach the conclusion of the theorem.
\end{proof}

\subsection[Traces on $\Cl^{a}_{\odd}(M)$ for $a\leq0$]{Traces on $\boldsymbol{\Cl^{a}_{\odd}(M)}$ for $\boldsymbol{a\leq0}$}\label{section:classificationZero}

In this section we assume as before, that the dimension $n$ is odd, and we prove that any trace on the algebra of odd-class operators of non-positive order is a linear combination of a generalized leading symbol trace and the canonical trace.

We can adapt Lemma~4.5 in \cite{leschneira} (see also Lemma 5.1.1 in~\cite{neirathesis}) in the case of odd-class operators:
\begin{lemma}\label{commCl02aaaodd}
If $a\in\Z$ is non-positive, then there exists an inclusion map
\[
[\Cl^{0}_{\odd}(M),\Cl^{2a}_{\odd}(M)]\hookrightarrow [\Cl^{a}_{\odd}(M),\Cl^{a}_{\odd}(M)],
\]
meaning that any commutator in $[\Cl^{0}_{\odd}(M),\Cl^{2a}_{\odd}(M)]$ can be written as a sum of commutators in $[\Cl^{a}_{\odd}(M),\Cl^{a}_{\odd}(M)]$.
\end{lemma}

\begin{proof}
By Lemma \ref{lema:oddClassproduct}, integer powers of an invertible dif\/ferential operator are odd-class ope\-ra\-tors. Hence we proceed as in the proof of Lemma~4.5 in \cite{leschneira} as follows:
Let $A\in \Cl^{0}_{\odd}(M)$, $B\in \Cl^{2a}_{\odd}(M)$. Consider a f\/irst-order positive def\/inite elliptic dif\/ferential operator $\Lambda$. For any $a\in\R$, $\Lambda^{a}$ and $\Lambda^{-a}$ are operators of order $a$ and $-a$, respectively, and therefore $A\Lambda^{a}$, $\Lambda^{a}A$, $\Lambda^{a}$, $B\Lambda^{-a}$,  $\Lambda^{-a}B$,  $AB\Lambda^{-a}$, $\Lambda^{-a}BA$ are operators in $\Cl^{a}_{\odd}(M)$. Moreover,
\begin{gather}
 [A\Lambda^{a},\Lambda^{-a}B]=AB-\Lambda^{-a}BA\Lambda^{a},\label{com1AB}\\
 [\Lambda^{a}A,B\Lambda^{-a}]=\Lambda^{a}AB\Lambda^{-a}-BA,\label{com2AB}\\
 [AB\Lambda^{-a},\Lambda^{a}]=AB-\Lambda^{a}AB\Lambda^{-a},\label{com3AB}\\
 [\Lambda^{-a}BA,\Lambda^{a}]=\Lambda^{-a}BA\Lambda^{a}-BA.\label{com4AB}
\end{gather}
Adding up the expressions in \eqref{com1AB}--\eqref{com4AB} yields twice the commutator $[A,B]$, so that the resulting expression belongs to the space of commutators $[Cl^{a}_{\odd}(M),Cl^{a}_{\odd}(M)]$.
\end{proof}

As in \eqref{ses1}, for a non-positive integer $a$ we also denote by $\pi_a$ the projection
\[
0\to \Cl^{2a-1}_{\odd}(M)\to \Cl^{a}_{\odd}(M)\overset{\pi_a}{\to}\Cl^{a}_{\odd}(M)/\Cl^{2a-1}_{\odd}(M)\to0,
\]
with corresponding splitting $\theta_a:\Cl^{a}_{\odd}(M)/\Cl^{2a-1}_{\odd}(M)\to \Cl^{a}_{\odd}(M)$, so that for any $A\in \Cl^{a}_{\odd}(M)$, $A-\theta_a(\pi_a(A))\in \Cl^{2a-1}_{\odd}(M)$.

The following result adds to the known classif\/ication of traces on pseudodif\/ferential ope\-ra\-tors, the classif\/ication of all traces on odd-class operators with f\/ixed non-positive order in odd-dimensions.

We f\/ix a non-positive integer~$a$, and describe any trace on $\Cl^{a}_{\odd}(M)$ (see Section~5.1.4 in~\cite{neirathesis}).

\begin{theorem}\label{classifictracesClodda}
If $a\in\Z$ is non-positive, any trace on $\Cl^{a}_{\odd}(M)$ can be written as a linear combination of a generalized leading symbol trace and the canonical trace.
\end{theorem}

\begin{proof}
Let $A\in \Cl^{a}_{\odd}(M)$. As in Corollary~\ref{corol2}, by Proposition~\ref{prop:PDOsumComutators} applied to $A-\theta_a(\pi_a(A))\in \Cl^{2a-1}_{\odd}(M)$, there exist operators $B_i\in \Cl^{0}_{\odd}(M)$, $C_i\in \Cl^{2a}_{\odd}(M)$, and a smoothing operator $R$ such that
\[
A-\theta_a(\pi_a(A))=\sum\limits_{i=1}^n[B_i,C_i]+R.
\]
By Lemma \ref{commCl02aaaodd}, there exist operators $D_1,\ldots,D_N,E_1,\ldots,E_N\in \Cl^{a}_{\odd}(M)$, such that
\begin{gather}
A-\theta_a(\pi_a(A))=\sum\limits_{k=1}^N[D_k,E_k]+R. \label{Cloddsumofcomm2}
\end{gather}
Applying $\TR$ to both sides of \eqref{Cloddsumofcomm2} yields
\[
\TR(A-\theta_a(\pi_a(A)))=\sum\limits_{k=1}^N\TR([D_k,E_k])+\TR(R)=\Tr_{L^2}(R).
\]
Hence, as in Theorem \ref{smooth}, for any pseudodif\/ferential idempotent $J$, of rank~1, there exist smoothing operators $S_1,\ldots,S_{N'}$, $T_1,\ldots,T_{N'}$, such that~(\ref{Cloddsumofcomm2}) becomes
\begin{gather}
A-\theta_a(\pi_a(A))=\sum_{k=1}^N[D_k,E_k] + \TR(A-\theta_a(\pi_a(A)))J+\sum_{j=1}^{N'}[S_j,T_j]. \label{comm2Clodd}
\end{gather}
Let $\tau:\Cl^{a}_{\odd}(M)\to\C$ be a trace on $\Cl^{a}_{\odd}(M)$. If we apply $\tau$ to both sides of~\eqref{comm2Clodd} we get
\begin{gather*}
\tau(A) =\tau(\theta_a(\pi_a(A)))+\TR(A-\theta_a(\pi_a(A)))\tau(J) \\
\phantom{\tau(A)}{}  =\tau(\theta_a(\pi_a(A)))-\TR(\theta_a(\pi_a(A)))\tau(J)+\TR(A)\tau(J).
\end{gather*}
So we conclude that $\tau$ is a linear combination of a generalized leading symbol trace and the canonical trace.
\end{proof}

\section{Determinants and traces}\label{section4}

We use the classif\/ication of traces on algebras of odd-class operators given in Theorem \ref{classifictracesClodda} to classify the associated determinants on the corresponding Fr\'echet--Lie group. Well-known ge\-ne\-ral results in the f\/inite-dimensional context concerning determinants associated with traces generalize to the context of Banach spaces (see~\cite{HarpeSkandalis}) and further to Fr\'echet spaces (see~\cite{lescure}).
\begin{definition}[Def\/inition 36.8 in~\cite{KrieglMichor}] \label{defin:ExponintialMapping}
A (possibly inf\/inite-dimensional) Lie group ${\mathcal{G}}$ with Lie algebra $\Lie({\mathcal{G}})$ admits an \emph{exponential mapping} if there exists a smooth mapping ${\Exp}:\Lie({\mathcal{G}})\to{\mathcal{G}}$ such that $t\mapsto {\Exp}(tX)$ is a one-parameter subgroup, i.e.\ a Lie group homomorphism $(\mathbb{R},+)\to {\mathcal{G}}$ with tangent vector $X$ at $0$.
\end{definition}

The existence of a smooth exponential mapping for a Lie group is ensured by a notion of regularity \cite{KrieglMichor,Milnor} on the group. Following \cite{KrieglMichor}, for J.~Milnor \cite{Milnor}, a Lie group ${\mathcal{G}}$ modelled on a~locally convex space is a regular Lie group if for each smooth curve $u:[0, 1]\to \Lie({\mathcal{G}}),$ there exists a smooth curve $\gamma_u:[0, 1]\to {\mathcal{G}}$ (which is unique: Lemma~38.3 in \cite{KrieglMichor}) which solves the initial value problem $\dot\gamma_u=\gamma_u\,u$ with $\gamma_u(0)=1_{\mathcal{G}},$ where $1_{\mathcal{G}}$ is the identity of ${\mathcal{G}},$ with smooth evolution map
\begin{gather*}
 C^{\infty}([0, 1], \Lie({\mathcal{G}})) \to {\mathcal{G}}, \qquad
 u \mapsto \gamma_u(1).
\end{gather*}
For example, Banach Lie groups (in particular f\/inite-dimensional Lie groups) are regular. If $E$ is a Banach space, then the Banach Lie group of all bounded
automorphisms of $E$ is equipped with an exponential mapping given by the series
\[
\Exp(X)=\sum_{i=0}^\infty\dfrac{X^i}{i!}.
\]

In \cite{kmoyIV} a wider concept of inf\/inite-dimensional Lie groups called regular Fr\'echet--Lie groups is introduced. In this paper we will consider the regular Fr\'echet--Lie groups of classical $\psi$DOs of non-positive order (see \cite{kmoy}).

If $\mathcal{G}$ admits an exponential mapping $\Exp$ and if a suitable inverse function theorem is appli\-cable, then $\Exp$ yields a dif\/feomorphism from a neighborhood of 0 in $\Lie({\mathcal{G}})$ onto a neighborhood of $1_{\mathcal{G}}$ in ${\mathcal{G}}$, whose inverse is denoted by $\Log$.

For our purpose in this section, we assume that the Lie group ${\mathcal{G}}$ is regular.

\begin{definition}[see Def\/inition~2 in \cite{lescure}]\label{DeftildeGdettrace}
Let $\mathcal{G}$ be a Lie group and let $\mathcal{\tilde{G}}$ be its subgroup of elements pathwise connected to the identity of $\mathcal{G}.$ A \emph{determinant} on $\mathcal{G}$ is a group morphism
\[
\Lambda: \ \mathcal{\tilde{G}}\to\mathbb{C}^*,
\]
i.e.\ for any $g, h \in \mathcal{\tilde{G}},$
\[
\Lambda(gh)=\Lambda(g)\Lambda(h).
\]
We also say that $\Lambda$ is multiplicative.

A \emph{trace} on the Lie algebra $\Lie(\mathcal{G})$ is a linear map $\lambda:\Lie(\mathcal{G})\to\C$, such that for all $u,v\in\mathcal{G}$,
\[
\lambda([u,v])=0.
\]
In our examples below $[u,v]=uv-vu$.
\end{definition}

The following lemma gives the construction of a locally def\/ined determinant on $\mathcal{G}$ from a~trace on $\Lie({\mathcal{G}})$.

\begin{lemma}[see Proposition 2 and Theorem~3 in \cite{lescure} which is based on \cite{HarpeSkandalis}]\label{lem:traceTOdeterminant}
A continuous trace $\lambda:\Lie({\mathcal{G}})\to\mathbb{C}$ gives rise to a determinant map $\Lambda:{\Exp}(\Lie({\mathcal{G}}))\subset{\mathcal{\tilde{G}}}\to\mathbb{C}^*$ def\/ined on the range of the exponential mapping by
\[
\Lambda(g):={\exp}(\lambda({\Log}(g))),
\]
where locally ${\Log}={\Exp}^{-1}$, making the following diagram commutative, for any small enough neighborhood $U_0$ of zero in $\Lie({\mathcal{G}})$:
\[
\begindc{\commdiag}[35]
\obj(1,3)[objLie(cal{G})]{$\hspace{1,5cm}U_0\subset \Lie({\mathcal{G}})$}
\obj(4,3)[objMathbb{C}]{$\mathbb{C}$}
\obj(4,1)[objMathbb{C}^*]{$\mathbb{C}^*$}
\obj(1,1)[objCal{G}]{${\Exp}(U_0)\subset {\mathcal{\tilde{G}}}$}
\mor{objLie(cal{G})}{objMathbb{C}}{$\lambda$}
\mor{objMathbb{C}}{objMathbb{C}^*}{${\exp}$}
\mor{objCal{G}}{objMathbb{C}^*}{$\Lambda$}
\mor{objLie(cal{G})}{objCal{G}}{${\Exp}$} [\atright,\solidarrow]
\enddc
\]
Moreover, $\Lambda$ is differentiable $($hence of class $C^1)$ at $1_{\mathcal{G}}$, with differential $D_{1_{\mathcal{G}}}\Lambda=\lambda$.
\end{lemma}

\begin{proof}
We f\/irst observe that for all $g\in{\Exp}(U_0)$, ${\log}(\Lambda(g))=\lambda({\Log}(g))$. Let $u\in U_0\subset \Lie({\mathcal{G}})$ be such that $g={\Exp}(u).$ Since ${\mathcal{G}}$ is a regular Lie group, we can consider the $C^1$-path $\gamma(t)={\Exp}(tu)$ going from $1_{\mathcal{G}}$ to ${\Exp}(u)=g.$ We have $\gamma(t)^{-1}\dot\gamma(t)=u$ and hence
\begin{gather}
\int_0^1\lambda\left(\gamma(t)^{-1}\dot\gamma(t)\right)dt=\lambda\left(\int_0^1\gamma(t)^{-1}\dot\gamma(t)dt\right)=\lambda(u)=\lambda({\Log}(g)) \label{eq:integralgammalog}
\end{gather}
using the continuity of $\lambda$ and that $\Log(1_{\mathcal{G}})=0$. It follows that if $\gamma_1$, $\gamma_2$ are two $C^1$-paths going from $1_{{\mathcal{G}}}$ to $g_1$ and $g_2$ respectively, then $\gamma_1\gamma_2$ is a $C^1$-path going from $1_{{\mathcal{G}}}$ to $g_1g_2$ and we have
\begin{gather*}
\lambda\big((\gamma_1(t)\gamma_2(t))^{-1}\dot{\widetilde{\gamma_1(t)\gamma_2(t)}}\big)=
\lambda\left(\gamma_2(t)^{-1}\gamma_1(t)^{-1}\dot{\gamma_1}(t)
\gamma_2(t) +\gamma_2(t)^{-1}\dot{\gamma_2}(t)\right)\\
\phantom{\lambda\big((\gamma_1(t)\gamma_2(t))^{-1}\dot{\widetilde{\gamma_1(t)\gamma_2(t)}}\big)}{}
= \lambda\left(\gamma_1(t)^{-1}\dot{\gamma_1}(t)\right)+\lambda\left(\gamma_2(t)^{-1}\dot{\gamma_2}(t)\right),
\end{gather*}
where we have used the tracial property of~$\lambda$.

 Now, for $g_1, g_2\in {\Exp}(U_0)\subset\mathcal{\tilde{G}}$,
\begin{gather*}
\log(\Lambda(g_1g_2)) =\lambda({\Log}(g_1g_2))
 =\lambda\left({\Log}(g_1)\right)+\lambda\left({\Log}(g_2)\right)
  =\log(\Lambda(g_1))+\log(\Lambda(g_2)),
\end{gather*}
and we can apply the map ${\exp}$ to both sides of this expression to reach the statement.
\end{proof}

Conversely, following~\cite{lescure} we give a construction of a trace from a determinant. Our proof here is dif\/ferent from the one in loc.\ cit.

\begin{lemma}[See Proposition~2 and Theorem~3 in \cite{lescure}]\label{lem:determinantTOtrace}
A determinant $\Lambda:\Exp(\Lie(\mathcal{G}))\to \C^*$, which is of class $C^1$ on ${\mathcal{G}}$, yields a continuous trace $\lambda:\Lie({\mathcal{G}})\to\mathbb{C}$ in the following way: for all $u\in \Lie({\mathcal{G}})$
\[
\lambda(u):=D_{1_{\mathcal{G}}}\Lambda(u)=\frac{d}{dt}\Big|_{t=0}\Lambda(\Exp(tu)),
\]
which makes the following diagram commutative:
\[
\begindc{\commdiag}[35]
\obj(1,3)[objCal{G}]{${\Exp(\Lie(\mathcal{G}))}$}
\obj(4,3)[objMathbb{C}^*]{$\mathbb{C}^*$}
\obj(4,1)[objMathbb{C}]{$\mathbb{C}$}
\obj(1,1)[objLie(cal{G})]{$\Lie({\mathcal{G}})$}
\mor{objLie(cal{G})}{objMathbb{C}}{$\lambda$}
\mor{objMathbb{C}}{objMathbb{C}^*}{${\exp}$}[\atright,\solidarrow]
\mor{objCal{G}}{objMathbb{C}^*}{$\Lambda$}
\mor{objLie(cal{G})}{objCal{G}}{${\Exp}$}
\enddc
\]
\end{lemma}

\begin{proof}
Let $u_1, u_2\in \Lie({\mathcal{G}})$. Then
\begin{gather*}
\lambda([u_1,u_2])
    = \frac{d}{ds}\Big|_{s=0}\, \frac{d}{dt}\Big|_{t=0}
       \Lambda\left({\Exp}(tu_1){\Exp}(su_2){\Exp}(-tu_1)\right)\\
\phantom{\lambda([u_1,u_2])}{}
=\frac{d}{ds}\Big|_{s=0}\, \frac{d}{dt}\Big|_{t=0}
       \Lambda({\Exp}(tu_1))\Lambda({\Exp}(su_2))\Lambda({\Exp}(-tu_1))\\
\phantom{\lambda([u_1,u_2])}{}
 = \frac{d}{ds}\Big|_{s=0}\, \frac{d}{dt}\Big|_{t=0}
       \Lambda({\Exp}(su_2)) =0.
\end{gather*}
Here we use the fact that $\Lambda$ is multiplicative, which implies that $\Lambda(g^{-1})=\Lambda(g)^{-1}$. The linearity of $\lambda$ can be proved using the commutativity of the diagram.
\end{proof}

\begin{remark}\label{rk:Tracedeterminant}
The two lemmata imply that continuous traces on $\Lie({\mathcal{G}})$ are in one to one correspondence with $C^1$-determinants on the open subset of ${\mathcal{G}}$ given by the range of the exponential mapping.
\end{remark}

In the following we assume that $M$ is an $n$-dimensional closed manifold and $n$ is odd. We are going to classify determinants on a neighborhood of the identity in the space of invertible pseudodif\/ferential operators $(I+\Cl^{a}_{\odd}(M))^*$, for a non-positive integer $a$. For that, let us f\/irst single out the Fr\'echet--Lie groups and Fr\'echet--Lie algebras we use.

The following proposition can be found in Proposition~3 in~\cite{lescure} for the case $(\Cl^{0}(M))^*$. See Proposition~6.1.4 in \cite{ouedraogothese} for the case of odd-class operators. We give here an exhaustive proof of this result since we did not f\/ind it in the literature.

\begin{proposition}\label{prop:LieGroupLieAlgebraCl^{0}_{odd}(M)}
$\left(\Cl^{0}_{\odd}(M)\right)^*$ is a Fr\'echet--Lie group which admits an exponential mapping and its Fr\'echet--Lie algebra is $\Cl^{0}_{\odd}(M)$.
\end{proposition}

\begin{proof}
By Lemma \ref{lema:oddClassproduct}, the composition of two operators in $\left(\Cl^{0}_{\odd}(M)\right)^*$ belongs to $\left(\Cl^{0}_{\odd}(M)\right)^*$ and the same holds for the inverse so that the set $\left(\Cl^{0}_{\odd}(M)\right)^*$ is a group in the Fr\'echet algebra $\Cl^{0}_{\odd}(M).$ Let us show that it is also an open subset of $\Cl^{0}_{\odd}(M)$. Let $A$ be an operator in $\left(\Cl^{0}_{\odd}(M)\right)^*$; we want to build an open neighborhood of $A$ in $\left(\Cl^{0}_{\odd}(M)\right)^*.$ The algebra $\Cl^{0}_{\odd}(M)$ is contained in $\Cl^{0}(M)$ which corresponds to all bounded classical $\psi$DOs on $L^2(M)$ and hence it is contained in the Banach algebra ${\mathcal{L}}(L^2(M))$ of bounded linear operators. By the local inverse function theorem, the set of invertible operators on ${\mathcal{L}}(L^2(M))$ is an open set. Hence $A$ admits an open neighborhood $V$ in the set of invertible operators in ${\mathcal{L}}(L^2(M))$.

 On the other hand (see~\cite{Shubin}), the inclusion $\Cl^{0}(M)\to {\mathcal{L}}(L^2(M))$ is continuous so that the inclusion $i:\Cl^{0}_{\odd}(M)\to\Cl^{0}(M)\to {\mathcal{L}}(L^2(M))$ is also continuous and the inverse image $i^{-1}(V)$ yields an open neighborhood of $A$ in $\left(\Cl^{0}_{\odd}(M)\right)^*.$ It follows that $\left(\Cl^{0}_{\odd}(M)\right)^*$ is canonically equipped with a manifold structure which makes it a Lie group.

  Let us now construct an exponential mapping on $\Cl^{0}_{\odd}(M)$. Given any operator $B$ in $\Cl^{0}_{\odd}(M)$, the dif\/ferential equation
\[
A^{-1}_t{\dot A}_t=B,\qquad A_0=I
\]
has a unique solution in $\left(\Cl^{0}(M)\right)^*$ given by
\[
A_t=\frac{i}{2\pi}\int_\Gamma {\exp}(t\lambda)(B-\lambda)^{-1}\,d\lambda,
\]
where $\Gamma$ is a contour around the spectrum of $B$. Note that $A_t$ is bounded since $B$ has zero order. Let us check that $A_t$ belongs to $\left(\Cl^{0}_{\odd}(M)\right)^*.$ The homogeneous components of the symbol of~$A_t$ are
\[
\sigma(A_t)_{-j}=\frac{i}{2\pi}\int_\Gamma {\exp}(t\lambda)b_{-j}(\lambda)\,d\lambda,
\]
where $b_{-j}(\lambda)$ denote the components of the resolvent $(B-\lambda)^{-1}$ of $B$ at the point $\lambda$. Explicitly we have
\begin{gather*}
b_{0}(\lambda) := (\sigma_0(B)-\lambda)^{-1},\qquad
b_{-j}(\lambda) := -b_{0}(\lambda)\sum_{\substack{l<j\\ k+l+\vert\alpha\vert=j}}\frac{(-i)^{|\alpha|}}{\alpha}\partial_{\xi}^\alpha\sigma_{-k}(B)\partial_{x}^\alpha b_{-l}(\lambda).
\end{gather*}
So that $(B-\lambda)^{-1}$ lies in $\Cl^{0}_{\odd}(M)$ since $B$ lies in $\Cl^{0}_{\odd}(M)$. It follows that $A_t$ lies in $\left(\Cl^{0}_{\odd}(M)\right)^*.$ This def\/ines an exponential mapping
\[
{\Exp}: \ \Cl^{0}_{\odd}(M) \to \left(\Cl^{0}_{\odd}(M)\right)^*.
\]
Moreover, it follows that for any smooth curve $u:[0, 1]\to \Cl^{0}_{\odd}(M),$ there exists a unique smooth curve $\gamma_u:[0, 1]\to \left(\Cl^{0}_{\odd}(M)\right)^*$ def\/ined by the following diagram
\[
\begindc{\commdiag}[35]
\obj(8,1)[objCl^{0}_{odd}(M)^*]{$\left(\Cl^{0}_{\odd}(M)\right)^*$}[\west]
\obj(4,1)[objCl^{0}_{odd}(M)]{$\Cl^{0}_{\odd}(M)$}[\west]
\obj(1,1)[obj01]{$[0, 1]$}[\west]
\mor{objCl^{0}_{odd}(M)}{objCl^{0}_{odd}(M)^*}{${\Exp}$}
\mor{obj01}{objCl^{0}_{odd}(M)}{$u$}
\enddc
\]
which solves the initial value problem $\gamma_u^{-1}\dot\gamma_u=u$, $\gamma_u(0)=I$.
\end{proof}

Inspired by Corollary 5.12 in \cite{kmoyIV}, and \cite{kmoy} we have the following

\begin{proposition}\label{frechetliegroup}
If $a<0$, the space of invertible odd-class $\psi$DOs
\[
\mathcal{G}:=(I+\Cl^{a}_{\odd}(M))^{*}=(\{I+A:A\in \Cl^{a}_{\odd}(M)\})^{*}
\]
is a Fr\'{e}chet--Lie group with Fr\'{e}chet--Lie algebra $\Cl^{a}_{\odd}(M)$, which admits an exponential mapping from $\Cl^{a}_{\odd}(M)$ to $(I+\Cl^{a}_{\odd}(M))^{*}$.
\end{proposition}

Explicitly this exponential mapping is given by
\begin{gather*}
\Exp: \ \Cl^{a}_{\odd}(M) \to\widetilde{\mathcal{G}},\qquad
A \mapsto\Exp(A)=\sum_{k=0}^\infty\dfrac{1}{k!}A^k.
\end{gather*}
This map restricts to a dif\/feomorphism from some neighborhood of the identity in $\Cl^{a}_{\odd}(M)$ to a neighborhood of the identity in $\widetilde{\mathcal{G}}$. The inverse is given by
\begin{gather*}
\Log: \ \widetilde{\mathcal{G}} \to \Cl^{a}_{\odd}(M),\qquad
I+A \mapsto\Log(I+A)=\sum_{k=1}^\infty\dfrac{(-1)^{k+1}}{k}A^k.
\end{gather*}

\subsection[Classification of determinants on $(I+\Cl^{a}_{\odd}(M))^*$ for $a\leq0$]{Classif\/ication of determinants on $\boldsymbol{(I+\Cl^{a}_{\odd}(M))^*}$ for $\boldsymbol{a\leq0}$}

As before we consider a non-positive integer number $a$. From the classif\/ication of traces on $\Cl^{a}_{\odd}(M)$ derived in Theorem~\ref{classifictracesClodda}, we infer a description of the determinants def\/ined on the range of the exponential mapping in $\left(I+\Cl^{a}_{\odd}(M)\right)^*$.

The following theorem supplements the known classif\/ication of determinants on classical pseudodif\/ferential operators.

\begin{theorem}\label{prop:ClassificationDetCl_{odd}^{0}(M)}
Let $a\in\Z$ be such that $a\leq0$. Determinant maps on $(I+\Cl^{a}_{\odd}(M))^*$ are given by a~two-parameter family: for any $c_1,c_2\!\in\!\C$, and for any linear map $\rho:Cl^{a}_{\odd}(M)/Cl^{2a-1}_{\odd}(M)\!\to\!\C$,
\[
\Det_{c_1,c_2}(\cdot)=\exp\big(c_1\rho\circ\pi_a(\Log(\cdot)) + c_2\TR(\Log(\cdot))\big).
\]
\end{theorem}

\begin{proof}
By Theorem \ref{classifictracesClodda}, any trace $\tau$ on $\Cl^{a}_{\odd}(M)$ is a linear combination of the canonical trace and a generalized leading symbol trace:
\[
\tau(\cdot)=c_1\rho\circ\pi_a(\cdot) + c_2\TR(\cdot),
\]
for some $c_1,c_2\in\C$, and for some linear map $\rho:Cl^{a}_{\odd}(M)/Cl^{2a-1}_{\odd}(M)\to\C$. Moreover, $\tau$ is continuous for the Fr\'echet topology of $\Cl^{a}_{\odd}(M)$. Then, applying Lemma \ref{lem:traceTOdeterminant} to ${\mathcal{G}}=\left(I+\Cl^{a}_{\odd}(M)\right)^*$ and $\Lie({\mathcal{G}})=\Cl^{a}_{\odd}(M)$, it follows that a determinant map on ${\mathcal{G}}$ is of the form
\begin{gather*}
\exp\big(c_1\rho\circ\pi_a(\Log(\cdot)) + c_2\TR(\Log(\cdot))\big).\tag*{\qed}
\end{gather*}
\renewcommand{\qed}{}
\end{proof}

The determinants given in Theorem \ref{prop:ClassificationDetCl_{odd}^{0}(M)} dif\/fer from the ones sometimes used by physicists for operators of the type $I+$\textit{Schatten-class\ operator} \cite{Mic,simon} which in contrast to these are not multiplicative but do extend the ordinary determinant (Fredholm determinant) for operators of the type $I+$\textit{trace-class operator}.

Here are some relevant specif\/ic cases
\begin{itemize}\itemsep=0pt
\item $\Det_{1,0}(\cdot)=\exp (\rho\circ\pi_a(\Log(\cdot)) )$ are extensions of the leading symbol determinants (see~\cite{rosenberg} for the case $\Cl^{0}(M)$).
\item $\Det_{0,1}(\cdot)=\exp (\TR(\Log(\cdot)) )$ is an extension of the Fredholm determinant (see Lemma~2.1 in \cite{scottres}, and \cite{simon}).
\end{itemize}

\subsection{Extension of determinants}\label{subsectionanotherconstructiondet}

Now we consider determinants constructed from a trace as above, by def\/ining both sides of Equation \eqref{eq:integralgammalog}, not only on a neighborhood of the identity in the range of the exponential mapping, but also on a set of admissible operators and on the pathwise connected component of the identity.

\begin{remark}
In this section, the word ``extension'' is in the sense that the determinants can be def\/ined on operators which not necessarily lie in the range of the exponential mapping. See Remark \ref{noteextension} below for the other sense of this word.
\end{remark}

\subsubsection{First extension}\label{subsectionfirstextension}

The f\/irst way to def\/ine these determinants is carried out by considering the right hand side of~\eqref{eq:integralgammalog}. In the case of ${\mathcal{G}}=\left(\Cl^{0}_{\odd}(M)\right)^*$ the logarithm can be def\/ined provided one chooses a~spectral cut $\theta$ thereby to f\/ix a determination $\log_\theta$ of the logarithm. We set
\begin{gather}
\Det^\lambda_\theta(A):= {\exp}\left(\lambda(\log_\theta A)\right). \label{Detlambda}
\end{gather}
Recall that if the operator $A$ lies in the odd-class and has even order, then the logarithm $\log_\theta A$ lies also in the odd-class. If $\phi$ is another spectral cut of $A$ such that $0\leq\theta<\phi<2\pi$, by formula~(\ref{eq:difflogarithms})
we have
\[
\log_\theta A-\log_\phi A= -2i\pi P_{\theta, \phi}(A),
\]
where $P_{\theta, \phi}(A)$ is an odd-class projection as in \eqref{projection}.

\begin{proposition}[Proposition~6.2.3 in~\cite{ouedraogothese}]\label{prop:FundGroupEntier}
Let $\lambda$ be any continuous trace on $\Cl_{\odd}^{0}(M)$. Suppose that $\lambda$ takes integer values on the image of $P_{\theta, \phi}(A)$ for all $\theta$, $\phi$ and $A$, where $A$ is an admissible operator with spectral cuts $\theta$ and $\phi$ as in~\eqref{eq:difflogarithms}. Then $\lambda$ gives rise to the map ${\Det}_\theta^\lambda$ in~\eqref{Detlambda}, on admissible operators, independent of the choice of the spectral cut~$\theta$.
\end{proposition}

Let us consider this construction for the traces given in Section~\ref{sectiontraces}: Let $A$ be an admissible operator in $\Cl_{\odd}^{0}(M)$  with spectral cut $\theta$.

 With the notation of Subsection \ref{subsectiongeneralizedleadingsymboltrace}, the determinant associated to the leading symbol trace $\lambda=\rho\circ\pi_0$ is def\/ined by
\[
 {\Det}_\theta^{\lambda}(A):={\exp}\left(\rho\circ\pi_0(\log_\theta A)\right).
\]
In Example 2 of \cite{lescure}, it is shown that if $P$ is a zero-order pseudodif\/ferential idempotent, then its leading symbol $p$ is also an idempotent so that the f\/ibrewise trace ${\tr}_x(p(x,\cdot))$ is the rank $\rk(p(x,\cdot)).$ Hence
\begin{gather*}
(\rho\circ\pi_0)(P_{\theta, \phi}(A)) =(\rho\circ\pi_0)\big({\tr}_x\big(\sigma^L(P_{\theta,\phi}(A)(x,\xi))\big)\big)
 =\rk\big(\sigma^L(P_{\theta, \phi}(A))\big)(\rho\circ\pi_0)(I).
\end{gather*}
It follows that ${\Det}_\theta^{\lambda}(A)$ is independent of the choice of the spectral cut $\theta$ for any linear map $\rho:Cl^{0}_{\odd}(M)/Cl^{-1}_{\odd}(M)\to\Z$.

 Observe that if $A$ is a zero-order odd-class operator so is $\log_\theta A$. Hence, the canonical trace extends to logarithms of admissible odd-class operators of order zero with its property of cyclicity in odd dimensions. The determinant associated to the canonical trace is def\/ined by
\begin{gather}
{\Det}_{\theta}^{\TR}(A):= {\exp}({{\TR}(\log_\theta A)}). \label{DetthetaTR}
\end{gather}
In contrast to the leading symbol trace, the canonical trace does not satisfy the requirement of Proposition~\ref{prop:FundGroupEntier} so that the associated determinant depends on the choice of spectral cut.

\subsubsection{Second extension}\label{subsectionsecondextension}

An alternative way to def\/ine these determinants is by considering the left hand side of~\eqref{eq:integralgammalog} and the expression for a determinant given in the proof of Lemma \ref{lem:traceTOdeterminant}:
\[
\Lambda(g)={\exp}\left(\int_0^1\lambda(\gamma(t)^{-1}\dot{\gamma}(t))\,dt\right),
\]
where $\gamma:[0,1]\to {\mathcal{G}}$ is a $C^1$-path with $\gamma(0)=1_{\mathcal{G}}$ and $\gamma(1)=g$.

Such an approach was adopted in \cite{HarpeSkandalis} by P.~de la Harpe and G.~Skandalis in  the case of a~Banach Lie group. In her thesis \cite{ducourtiouxthese}, C.~Ducourtioux adopted this point of view to construct a determinant associated to a weighted trace with associated Lie algebras $\Cl^{0}(M)$ and $\Cl(M)$. In \cite{lescure}, J.-M.~Lescure and S.~Paycha showed that such a construction extends to Fr\'echet--Lie groups with exponential mapping.

As in Def\/inition \ref{DeftildeGdettrace}, let $\widetilde{\mathcal{G}}$ denote the pathwise connected component of the identity $1_{\mathcal{G}}$ of~${\mathcal{G}}$ and~${\mathcal{P}}(\mathcal{G})$ the set of $C^1$-paths $\gamma:[0,1]\to{\mathcal{G}}$ starting at $1_{\mathcal{G}}$ ($\gamma(0)=1_{\mathcal{G}}$) in $\widetilde{\mathcal{G}}$. On ${\mathcal{P}}(\mathcal{G})$ we introduce the map: $\Det_\lambda:{\mathcal{P}}(\mathcal{G})\to \mathbb{C}^*$ def\/ined by{\samepage
\begin{gather}
\Det^\lambda(\gamma):= {\exp}\left( \int_\gamma \lambda(\omega)\right)={\exp}\left(\int_0^1 \lambda(\gamma^*\omega)\right), \label{Detlambdagamma}
\end{gather}
where $\omega= g^{-1} dg$ is the Maurer--Cartan form on ${\mathcal{G}}$.}

Since $\lambda$ satisf\/ies the tracial property, $\Det^\lambda$ has the multiplicative property:

\begin{lemma}\label{lema:MultiplicativityDet}
Let $\gamma_1$, $\gamma_2$ be two $C^1$-paths in ${\mathcal{P}}(\mathcal{G})$. Then
\[
\Det^\lambda(\gamma_1\gamma_2)=\Det^\lambda(\gamma_1)\Det^\lambda(\gamma_2).
\]
\end{lemma}

\begin{proof}
The same proof applies as in Lemma \ref{lem:traceTOdeterminant}.
\end{proof}

In general the Maurer--Cartan form  $\omega= g^{-1} dg$ is not  exact on ${\mathcal{G}}$ so that for a $C^1$-path $c:[0,1] \to {\mathcal{G}}$ with $c(0)=c(1),$ the integral $\int_c\omega=\int_0^1c^*\omega$ does not vanish.

\begin{proposition}[Proposition~6.2.2 in~\cite{ouedraogothese}]\label{prop:detpath}\qquad
\begin{enumerate}\itemsep=0pt
\item[$1.$] Let ${\mathcal{P}}_{\rm cl}({\mathcal{G}})$ denote the space of closed $C^1$-paths (loops) in ${\mathcal{G}}$. The map
\begin{gather*}
\Phi: \ {\mathcal{P}}_{\rm cl}({\mathcal{G}}) \to \Lie({\mathcal{G}}), \qquad
c \mapsto \int_c\omega=\int_0^1c^*\omega
\end{gather*}
induces a map $\Phi: \Pi_1({\mathcal{G}})\to\Lie({\mathcal{G}})$ on the fundamental group $\Pi_1({\mathcal{G}})$ of ${\mathcal{G}}$. Consequently, the map $\Det^\lambda$ defined in~\eqref{Detlambdagamma} only depends on the homotopy class of the path $\gamma$.
\item[$2.$] If $\Det^\lambda(\Pi_1({\mathcal{G}}))=1$, then it induces a determinant map:
\begin{gather*}
 \Det^\lambda:\ \widetilde{\mathcal{G}} \to \C^*, \qquad
 g \mapsto {\exp}\left(\int_0^1 \lambda(\gamma^*\omega)\right)
\end{gather*}
independently of the choice of path $\gamma$.
\item[$3.$]
 If $g$ lies in the range of the exponential mapping ${\Exp}$ then
\[
\Det^\lambda(g)= {\exp}\left(\lambda({\Log}(g))\right),
\]
where $ {\Log}= {\Exp}^{-1}$ is the inverse of the exponential mapping.
\end{enumerate}
\end{proposition}

\begin{remark}\label{noteextension}
Item 3 of this proposition shows that $\Det^\lambda$ is an extension of the determinant map def\/ined in Lemma~\ref{lem:traceTOdeterminant}.
\end{remark}

\begin{proof}
1.~We want to show that two homotopic loops $c_1$ and $c_2$ have common primitive. Let us f\/irst recall the following general construction of a  primitive: for $\omega$ a dif\/ferential form on ${\mathcal{G}}$, let  $\gamma:[0,1]\to {\mathcal{G}}$ be a $C^1$-path and $F:[0,1]\to {\mathcal{G}}$  be such that for any $t\in [0,1]$, $F^\prime(t)=\omega(\gamma(t))\gamma^\prime(t).$ If $\omega$ is an exact form i.e.\ $w=df$ for some $f\in\mathcal{G}$, then $F(t)=f(\gamma(t))$ is a primitive of $F^\prime.$

 If the form $\omega$ is closed, then $\omega$ is locally exact. In this case, let $0=t_0<t_1<\cdots<t_k=1$ be a~subdivision of the interval $[0,1]$ such that $\gamma([t_{i-1},t_i])$ is a subset of ${\mathcal{G}}$, and there exists $f_i$ def\/ined on $[t_{i-1},t_i]$ such that $df_i=\omega$. We can construct a function $F(t)$ on $[0,1]$ in the following way: $F(t)=f_0(\gamma(t))$ on $[t_0,t_1],$ $F(t)=f_1(\gamma(t))-h_1$ on $[t_{1},t_2]$ where $h_1=f_1(\gamma(t_1))-f_{0}(\gamma(t_1))$ and for $i=3,\dots,k$, $F(t)=f_i(\gamma(t))-h_i$ on $[t_{i-1},t_i]$ where $h_i=f_i(\gamma(t_i))-f_{i-1}(\gamma(t_i))+h_{i-1}$.

Now let $F(t)$ and $G(t)$ be primitives  of $c_1$ and $c_2$ respectively. Since $c_1$ and $c_2$ are homotopic, there exists a family of $C^1$-paths $(\alpha_i)_{0\leq i\leq k}$ def\/ined in a neighborhood of $1_{\mathcal{G}}$ such that $c_1=c_2\prod\limits_{i=0}^k\alpha_i$. Each path $\alpha_i$ is closed so that $\int_c\omega$  vanishes on $\alpha_i$. It follows that $F(t)=G(t)$, and therefore the map $\Phi$ is well-def\/ined on $\Pi_1({\mathcal{G}}).$

 2.~The multiplicativity of $\Det^\lambda$ on ${\mathcal{G}}$ follows from Lemma \ref{lema:MultiplicativityDet}: If $g_1$, $g_2$ are two elements of ${\mathcal{G}}$ and $\gamma_1$, $\gamma_2$ are two $C^1$-paths in ${\mathcal{P}}(\mathcal{G})$ such that $\gamma_1(1)= g_1$ and $\gamma_2(1)= g_2$, then $\Det^\lambda(g_1g_2)=\Det^\lambda(g_1)\Det^\lambda(g_2)$.

 3.~For $g$ in the range of the exponential mapping and $\gamma$ a $C^1$-path in ${\mathcal{P}}(\mathcal{G})$ such that $\gamma(1)= g$, ${\Log}= {\Exp}^{-1}$ is well-def\/ined so that $\lambda({\Log}(\gamma(t)))$ is a primitive of $\lambda(\gamma(t)^{-1}\dot{\gamma}(t))$. Thus we recover equation~\eqref{eq:integralgammalog}:
\begin{gather*}
\int_0^1\lambda(\gamma(t)^{-1}\dot{\gamma}(t))\,dt=\lambda({\Log}(g)).\tag*{\qed}
\end{gather*}
\renewcommand{\qed}{}
\end{proof}

From Proposition \ref{prop:FundGroupEntier} we can see that in the f\/irst extension of the def\/inition of a determinant, the non-dependency on the spectral cut is controlled by the image of the projection $P_{\theta, \phi}$ by the trace~$\lambda$. From Proposition \ref{prop:detpath}, the map $\Det^\lambda$ is well def\/ined if the image of the fundamental group of $\left(\Cl^{0}_{\odd}(M)\right)^*$ is trivial.
The following theorem shows that in both def\/initions (Subsections~\ref{subsectionfirstextension} and~\ref{subsectionsecondextension}) the non-dependency on the spectral cut and the multiplicativity of $\Det^\lambda$ rely on the image of the fundamental group of $\left(\Cl^{0}_{\odd}(M)\right)^*$. It provides a way to identify which of the maps $\Det^\lambda$ def\/ined above send this fundamental group to $1$ and therefore are determinants in the sense of Def\/inition~\ref{DeftildeGdettrace}.

{\sloppy
\begin{theorem}[see Section~4.5 in~\cite{kontsevichdeterminant} and Lemma~A.5 in \cite{ducourtiouxthese}]
The fundamental group of $\left(\Cl^{0}_{\odd}(M)\right)^*$ is generated by the homotopy classes of loops $\{{\Exp}(2i\pi tP)\}_{0\leq t\leq1}$, where $P$ is a projector in $\Cl^{0}_{\odd}(M)$.
\end{theorem}

}

\begin{remark}
Let us give some comments for the case $\Cl_{\odd}(M)$. As seen in Theorem~\ref{prop:UniquenessTrace}, any trace in $\Cl_{\odd}(M)$ is proportional to the canonical trace. For the subalgebra $\Cl_{\odd}^{0}(M)$ in~\eqref{DetthetaTR} we def\/ined the determinant associated to the canonical trace of an operator $A$ with spectral cut $\theta$. Unfortunately, this def\/inition cannot be extended a priori to positive order elements of $\Cl_{\odd}(M)$. Indeed, as we said before (see Subsection \ref{subsectionfirstextension}), if the order of $A$ is even, for any spectral cut $\theta$ of $A$, the logarithm $\log_\theta A$ is also odd-class, so the canonical trace extends to logarithms of odd-class operators with even order and one can extend the determinant by: ${\Det}_{\theta}^{\TR}(A):= {\exp}({{\TR}(\log_\theta A)})$. This is no longer true if the order of~$A$ is odd. M.~Braverman in~\cite{Braverman} introduced the notion of symmetrized trace to def\/ine a symmetrized determinant on odd-class operators. It was shown in~\cite{ouedraogo} that this symmetrized trace is in fact the canonical trace. The idea is to compute the average of the usual terms given by two spectral cuts of~$A$, namely~$\theta$ and~$\theta-a\pi$, where~$a$ is the order of~$A$; one obtains the following def\/inition of symmetrized logarithm:
\[
\log_\theta^{\rm sym}A:=\frac{1}{2} (\log_\theta A+\log_{\theta-a\pi}A ).
\]
This symmetrized logarithm is also odd-class and, once again, one can use the canonical trace to def\/ine a determinant by setting:
\[
{\rm DET}_{\theta}^{\rm sym}(A):= {\exp}\big({\TR}\big(\log_\theta^{\rm sym} A\big)\big).
\]
Notice that if $a=0$, ${\rm DET}_{\theta}^{\rm sym}(A)={\Det}_{\theta}^{\TR}(A)$.
However this symmetrized determinant also depends on the spectral cut $\theta$ and under suitable assumptions it is multiplicative up to a sign (see~\cite{Braverman}).
\end{remark}

\subsection*{Acknowledgements}

We are very grateful to Professor Sylvie Paycha for her advice and encouragement along an important part of our careers.  We express our gratitude to the Max Planck Institute for Mathe\-matics in Bonn and the University of Regensburg where part of this work was done.  We would also like to thank the anonymous referees for providing constructive comments on the manuscript.

\pdfbookmark[1]{References}{ref}
\LastPageEnding

\end{document}